\documentclass[11pt]{amsart}
\usepackage{amsmath,amsthm,amsfonts,amssymb,bbm, setspace,graphicx,float,subfigure,xcolor,mathrsfs,mathtools, Style}
\usepackage{todonotes}
\usepackage[T1]{fontenc}
\usepackage[utf8]{inputenc}
\usepackage{enumitem}
\setlist[enumerate,1]{label={(\roman*)}}

\newcommand{\conv}[1]{\ensuremath{\operatorname{conv}(W #1)}}

\begin{document}

\author{Lasse L. Wolf}
\address{Institut des Hautes Études Scientifiques, 35 Route de Chartres,
91440 Bures-sur-Yvette} \email{\href{mailto:wolf@ihes.fr}{wolf@ihes.fr}}

\newcommand{\revision}[1]{{\color{blue} #1}}
\newcommand{\Lasse}[1]{{\color{red} \sf $\clubsuit\clubsuit\clubsuit$ Lasse: [#1]}}

\title[Limit cone and bounds on growth function]{The Limit Cone and Bounds on the Growth Indicator Function}

\begin{abstract}
	Given a real semisimple Lie group $G$ with finite center and a discrete subgroup
$\Gamma \subset G$ whose limit cone is does not contain the extremal rays of the Weyl chamber
we show that Quint's growth indicator function $\psi_\Gamma$ is bounded by the half sum of positive roots $\rho$,
i.e.~it has slow growth, implying that the representation $L^2(\Gamma \backslash G)$ is tempered.
In particular, this holds for each $I$-Anosov subgroup
provided that $I$ contains at least two distinct simple roots
that are not interchanged by the opposition involution.
\end{abstract}

\subjclass[2020]{22E40, 22E46, 58C40}

\setcounter{tocdepth}{1} 

\maketitle

\section{Introduction}

Let $G$ be a real semisimple Lie group
with finite center
and $\Gamma\subset G$ a discrete subgroup.
In this paper we study the growth rate of $\Gamma$ in $G$ in the case
where the real rank of $G$ is at least $2$.
Since for lattices the growth rate of $\Gamma$ in $G$
is always comparable with the volume growth in $G$,
this question becomes interesting as soon as
$\Gamma$ has infinite covolume
as abundant previous works demonstrate \cite{pattersonlimitset,Ben97,Qui02,Sam}
and we will assume this from now on.

To define the precise notion encoding the growth rates,
we need to establish some notation.
Recall that $G$ admits a Cartan decomposition $G=K \exp(\mathfrak a_+) K$.
For $g\in G$ we denote by $\mu_+(g)$ the $\mathfrak{a}_+$-component of $g$
which can be thought of as a higher dimensional distance of $gK$ to the origin
$eK$ in the symmetric space $G/K$.
The critical exponent
\[
	\delta_\Gamma = \inf \left\{s \geq 0\colon 
	\sum_{\gamma\in \Gamma} e^{-s \|\mu_+(\gamma)\|}<\infty\right\}
\]
has been refined by Quint \cite{Qui02} in higher rank.
He introduced the growth indicator function $\psi_\Gamma\colon \mathfrak{a} \to \R\cup\{-\infty\}$:
\begin{align*}
\psi_\G(H)\coloneqq \|H\| \inf_{\mathcal C\ni H} \inf\left\{s\in\R\colon \sum_{\g\in \G,\mu_+(\g)\in\mathcal C} e^{-s\|\mu_+(\g)\|} <\infty\right \},
	\end{align*}
where the first infimum runs over all open cones $\mathcal C\subseteq\mathfrak{a}$ with $H\in \mathcal{C}$.
We recall some of its properties in Subsection \ref{ss:gif}.
The growth indicator function $\psi_\Gamma$ is closely related to the critical exponents
\[
	\delta_\mu= \inf\left \{s\in \R \colon 
	\sum_{\gamma\in \Gamma} e^{-s \mu ( \mu_+(\gamma))}<\infty\right\}
\]
for $\mu \in \mathfrak{a}^\ast$ (see Subsection~\ref{ss:critexp}).

Another object related to $\psi_\Gamma$ is the limit cone $\mathcal{L}_\Gamma$
which was introduced by Benoist \cite{Ben97}.
Its importance lies in its connection with proper actions of $\Gamma$ on homogeneous spaces,
see \cite{DK, KT, DO} for recent results.
$\mathcal{L}_\Gamma$ is the asymptotic cone of $\mu_+(\Gamma)$,
i.e.~
\[
	\mathcal{L}_\Gamma = \{v\in \mathfrak{a}_+\colon |\mathcal{C}\cap \mu_+(\Gamma)| =\infty \text{ for all open cones }\mathcal{C} \text{ containing } v\},
\]
and coincides with the cone
where $\psi_\Gamma$ is finite:
\[
	\mathcal{L}_\Gamma = \{v\in \mathfrak{a}_+ \colon \psi_\Gamma(v) > -\infty\}.
\]
Note that $\mathcal{L}_\Gamma$ only encodes the directions in the Weyl chamber $\mathfrak{a}_+$
in which infinitely many $\Gamma$-points occur
without any information on the growth rate. The latter is encoded in the values of the growth indicator function $\psi_\Gamma$
and the critical exponents.

In this article, we show that nonetheless the location of the limit cone
inside the positive Weyl chamber $\mathfrak{a}_+$
has strong implications for the size of $\psi_\Gamma$ and the critical exponents.

To make this more precise, let $\Pi$ be the set of simple (restricted) roots
associated to the Weyl chamber $\mathfrak{a}_+$
and $\iota$ the opposition involution.
Furthermore, 
let $\rho\in \mathfrak{a}^\ast$ be the half sum of positive roots
counted with multiplicities (see Subsection~\ref{sec:not}).
Recall that $\mathfrak{a}_+ = \{v \in \mathfrak{a}\colon \alpha(v)\geq 0\: \forall\alpha\in \Pi\}$.
Let us choose extremal vectors $v_\alpha$ of $\mathfrak{a}_+$,
i.e.~$\beta(v_\alpha)=0$ for $\beta\in \Pi\setminus\{\alpha\}$,
and we normalize them by $\|v_\alpha\| =1$.
Let $C_\alpha$ be the cone $\R_{\geq 0} v_\alpha + \R_{\geq 0} v_{\iota \alpha}$
spanned by $v_\alpha$ and $v_{\iota \alpha}$.
Our main theorem is the following.
\begin{theorem}
	\label{mainintro}
	If $C_ \alpha \cap \mathcal{L}_\Gamma = \{0\} $ for all $\alpha \in \Pi$,
then
\[
	\psi_\Gamma(v) \leq \rho(v)
	\quad\text{for all}\quad v\in \mathfrak{a}_+.
\]
\end{theorem}

Before coming to other results, let us comment on this one.
We first remark that, as written in Corollary~\ref{cor:twowallsavoided},
we can relax the assumption
on $\mathcal{L}_\Gamma$ by requiring the same for the open subcone
\[
	\mathcal{L}_\Gamma' \coloneqq \{v\in \mathcal{L}_\Gamma \colon \psi_\Gamma(v) >\rho(v)\}.
\]
This should be thought of the limit cone for the modified growth indicator function
\[
	\psi_\Gamma' \coloneqq \psi_\Gamma -\rho
\]
which appears naturally in our study.
In particular, in the case where the opposition $\iota$ is trivial, we have the equivalence
\begin{align*}
	\psi_\Gamma \leq \rho \iff \psi_\Gamma (v_\alpha) \leq \rho(v_\alpha) \quad \forall \alpha\in \Pi.
\end{align*}

Theorem~\ref{mainintro} can be applied to Anosov subgroups and
we shortly recall their definition.
For a non-empty subset $I \subseteq \Pi$,
a finitely generated subgroup $\Gamma$ is called \emph{$I$-Anosov}
if there exist constants $C,c$ such that, for all $\alpha\in I$,
\[
	\alpha(\mu_+(\gamma)) \geq C |\gamma| -c
	\quad \text{for all}\quad
	\gamma\in \Gamma,
\]
where $|\gamma|$ is the word length of $\gamma$
with respect to some fixed generating set.
A $\Pi$-Anosov subgroup is also called Borel Anosov.
It is clear that if $\Gamma$ is $I$-Anosov,
then it is ${I\cup \iota I}$-Anosov.
Furthermore, it is immediate from the definition
that, for $\alpha\in I$, $\mathcal{L}_\Gamma \cap \ker \alpha =\{0\}$
and therefore $C_\beta \cap \mathcal{L}_\Gamma =\{0\}$ for $\beta\notin \{\alpha,\iota \alpha\}$.
Hence, we get the following corollary.
\begin{corollary}
	\label{cor:Anosovintro}
	If $\Gamma$ is an $I$-Anosov subgroup with $|I/(\alpha\sim \iota \alpha)|\geq 2$,
	then $\psi_\Gamma \leq \rho$.
\end{corollary}
We note that 
the assumption of $\Gamma$ being $I$-Anosov is much stronger than the assumption on the limit cone
in Theorem~\ref{mainintro}.
There are many more examples satisfying the limit cone assumption
such as relatively Anosov subgroups \cite{CZZ22} which are not Anosov.
We also note that there are examples of $I$-Anosov subgroups
with $\psi_\Gamma \not \leq \rho$
where $I$ contains only one simple root, see Subsection~\ref{ss:nontempex}.
In this sense, Corollary~\ref{cor:Anosovintro} is optimal.

Let us now elaborate on the property $\psi_\Gamma \leq \rho$.
This inequality first appeared in \cite{OhTemperedness}
where the authors showed that it is equivalent to the temperedness\footnote{see Subsection~\ref{sec:tempered} for the definition.}
 of the $G$-representation
$L^2(\Gamma\bk G)$ for Borel Anosov subgroups using a local mixing
result \cite{ELO23}.
Recently, this equivalence has been shown to hold for any discrete subgroup in \cite{LWW}.
In fact, $\psi_\Gamma\leq \rho$
is also equivalent to spectral properties of the Laplacian
and the full algebra of $G$-invariant differential operators on $G/K$ acting
on the locally symmetric space $\Gamma\bk G/K$,
as well as dynamical properties of $G\curvearrowright \Gamma\bk G$ (see Corollary~1.3 therein).

These equivalent properties were previously known to hold only for Hitchin subgroups of $\mathrm{SL}_n(\R)$
by the work of
\cite{OhTemperedness}
using \cite{PS17}.
Kim, Minsky, and Oh \cite{KMO24} then showed that $\psi_\Gamma\leq \rho$
holds for any Borel Anosov subgroup if $G$ is a product of rank one group
(see also \cite{WW23})
which led to the conjecture that $\psi_\Gamma \leq \rho$
should hold for any Borel Anosov subgroup in a higher rank group $G$
\cite[Conjecture~1.14]{KMO24}.
This conjecture has been proven in \cite{LWW}
except for the Lie algebras $\mathfrak{sl}_3(\mathbb{K})$, $\mathbb{K} = \R,\C,\mathbb{H}$, and $\mathfrak{e}_6^{-26}$.
For non-Borel Anosov subgroups further progress was made in \cite{DKO}
where bounds on $\psi_\Gamma$ in terms of the limit set were obtained
yielding a criterion for $\psi_\Gamma \leq \rho$
for $I$-Anosov subgroups in terms of the dimension of the limit set.
Their criterion produced many new examples where the limit set is known
\cite{PSW23, CZZ22, BILW}.
Corollary~\ref{cor:Anosovintro} considerably strengthens these results
by requiring the Anosov property only for two roots,
which are not in the same $\iota$-orbit,
and without imposing any additional assumptions on the limit set.
Nevertheless,
the conjecture of Kim-Minsky-Oh on Borel Anosov subgroups remains open
for the Lie algebras listed above as one has $\Pi = \{\alpha,\iota\alpha\}$,
i.e.~$|\Pi/(\alpha\sim \iota\alpha)| =1$ in these cases.

\subsection{Further results}
\label{ss:furtherresults}
In order to state additional results that we obtain in this article,
let us for simplicity assume that $G$ is an algebraic group
and $\Gamma$ is Zariski dense in $G$.
In this case, $\mathcal{L}_\Gamma$ is convex, $\psi_\Gamma$ is concave,
and there is a unique unit vector $v_\Gamma\in \mathcal{L}_\Gamma$ of maximal growth,
i.e.~$\max_{\|v\|=1} \psi_\Gamma (v) = \psi_\Gamma(v_\Gamma)$.
The same holds for the modified growth indicator function $\psi_\Gamma' = \psi_\Gamma -\rho$
if $\psi_\Gamma \not \leq \rho$
(i.e.~$\psi_\Gamma(v) > \rho(v)$ for some $v\in \mathfrak{a}_+$,
see Proposition~\ref{prop:mu_is_biggest_growth}):
\[
	\max_{\|v\|=1} (\psi_\Gamma-\rho)(v) = (\psi_\Gamma - \rho)(v'_\Gamma)\eqqcolon \delta'
\]
for a unique unit vector $v_\Gamma'$.
We prove the following theorem.
\begin{theorem}
	\label{thm:oneweight}
	If $\psi_\Gamma \not \leq \rho$
	and $v_\alpha + v_{\iota \alpha} \not \in \mathcal{L}'_\Gamma$ for some $\alpha\in \Pi$,
	then $v_\Gamma' \in \ker \alpha\cap \ker \iota \alpha$.
\end{theorem}
We note that $v_\alpha + v_{\iota \alpha}\notin \mathcal{L}'_\Gamma$ is equivalent to
$C_\alpha \cap \mathcal{L}'_\Gamma = \emptyset$ as $\mathcal{L}_\Gamma'$ is convex.

This theorem directly implies Theorem~\ref{mainintro} by considering all simple roots:
If $C_\alpha \cap \mathcal{L}'_\Gamma = \emptyset$ for all $\alpha$, then
$
v'_\Gamma \in \bigcap_{\alpha\in \Pi} \ker \alpha = \{0\}
	$
	which
	contradicts $\|v_\Gamma'\| =1$.
	Hence $\psi_\Gamma \not\leq \rho$ cannot hold.

	By considering all but one simple root we get the following corollary.
	This applies for example to $\{\alpha\}$-Anosov subgroups.
\begin{corollary}
	\label{thm:onewallintro}
	Let $\alpha\in \Pi$.
	If $\psi_\Gamma \not \leq \rho$ and
	$v_\beta + v_{\iota\beta}\notin \mathcal{L}'_\Gamma$ for all $\beta\in \Pi\setminus \{\alpha, \iota \alpha\}$
	(e.g.~if $\mathcal{L}'_\Gamma$ is disjoint from the facet $\ker \alpha \cap \mathfrak{a}_+$),
then
\[
	v'_\Gamma = \frac{1}{\|v_\alpha + v_{\iota\alpha}\|} (v_\alpha+ v_{\iota\alpha}),
\]
i.e.~$v_\Gamma'$ is the unique unit vector in $\mathfrak{a}_+$
that is $\iota$-invariant and satisfies $\beta(v_\Gamma') = 0$ for all $\beta\in \Pi\setminus \{\alpha,\iota\alpha\}$.

In particular, if $\iota \alpha =\alpha$, then $v'_\Gamma = v_\alpha$.
\end{corollary}

For a subset $I\subseteq\Pi$ let $\mathfrak{a}_I \coloneqq \bigcap_{\alpha\in \Pi\setminus I} \ker \alpha = \langle v_\alpha \colon \alpha\in I\rangle$.
The following theorem strengthens the contraposition of Theorem~\ref{thm:oneweight}
by determining $\psi_\Gamma$ on $v_\alpha + v_{\iota \alpha}$.

\begin{theorem}
	\label{thm:psilinearintro}
	Let $I=\{\alpha \in \Pi \colon  \alpha(v'_\Gamma) >0\}$.
Then
\[
	\psi_\Gamma(v) = \rho(v) +\delta' \langle v'_\Gamma,v\rangle
\]
for all $v\in \mathfrak{a}_I \cap \{\iota =1\}$.
In particular, $\mathfrak{a}_I \cap \{\iota =1\}$ is contained in the modified limit cone $\mathcal{L}'_\Gamma$.
\end{theorem}

\subsection{Spectral theory}
As we describe in Section~\ref{sec:outline},
the proofs rely on spectral theoretic results.
As a byproduct of our analysis,
we obtain properties of the joint spectrum
(see Section~\ref{ssec:jointspec} for the definition)
that we describe now.
Let $\Delta$ be the (non-negative) Laplace-Beltrami operator on a locally symmetric space
$\Gamma \bk G/K$.
It is of interest to determine the spectrum of $\Delta$ as an operator on $L^2(\Gamma \bk G/K)$.
Important results in this area were obtained by
Elstrodt \cite{MR360472, MR360473, MR360474}, Patterson \cite{pattersonlimitset}, Corlette \cite{Cor90} in rank one,
and Leuzinger \cite{Leu04}, Weber \cite{Web08}, Anker-Zhang \cite{AZ22}, Weich together with the author \cite{WW23b}, the author with Zhang \cite{WZ23} in higher rank.
One main result is to connect the bottom of the Laplace spectrum $\min \sigma(\Delta)$
with the critical exponents \cite[Corollary~1.3]{WZ23}:
\begin{equation}
	\label{eq:bottomLaplace}
	\min \sigma(\Delta) = \|\rho\|^2 - \max(0,\delta')^2.
\end{equation}

In higher rank however, the role of the Laplacian is played by the algebra $\mathbb{D}(G/K)$
of $G$-invariant differential operators on $G/K$ acting on $\Gamma \bk G/K$.
One can define a joint spectrum $\wt\sigma_\Gamma \subseteq \mathfrak{a}_\C^\ast$ (see Section~\ref{sec:spherical_dual})
such that $\sigma(\Delta)$ is the closure of $\{\|\rho\|^2 -\|\Re \lambda\|^2+\|\Im \lambda\|^2\colon \lambda\in \wt \sigma_\Gamma\}$.
In rank one, \eqref{eq:bottomLaplace} transforms to $\delta' = \max_{\lambda\in \wt \sigma} \|\Re \lambda\|$.
In this case, $\Re \lambda=0$ or $\Im \lambda=0$ so that $\delta' \in \wt \sigma_\Gamma$
as $\wt \sigma_\Gamma$ is closed.
In higher rank, we get the same result.
\begin{theorem}
	\label{thm:spectralintro}
	If $\Gamma$ is Zariski dense, then
	\[
		 \langle\max(0,\delta')  v_\Gamma', \cdot\, \rangle \in \wt \sigma_\Gamma
	\]
	and 
	\[
		\Re \lambda(v) \leq \max(0,\delta') \langle v,v'_\Gamma\rangle
	\]
	for all $\lambda\in \wt \sigma_\Gamma$ and $v\in \mathfrak{a}_+$.
\end{theorem}

\subsection{Related results}

\subsubsection{Critical exponents}
\label{ss:critexp}
The study of critical exponents has received increasing attention in recent years
due to their various applications including
description of limit sets, definition of Patterson-Sullivan measures,
and formulation as entropy, see e.g. \cite{Qui02,DK,CZZ22,OMT23, DKO}.
Hence, gaining knowledge about the critical exponents in general
and for specific examples is of great importance.
In the proof of Theorem~\ref{mainintro}
we use a modified $\mu$-critical exponent
(for $\mu\in \mathfrak{a}^\ast$ which is positive on $\mathcal{L}_\Gamma$)
defined by
\[
	\delta'_\mu= \inf\left \{s\in \R \colon 
	\sum_{\gamma\in \Gamma} e^{-(s \mu+\rho) ( \mu_+(\gamma))}<\infty\right\}=
	\max_{v\in \mathfrak{a}_+} \frac{(\psi_\Gamma-\rho)(v)}{\mu(v)}.
\]
Hence, in the setting of Theorem~\ref{mainintro},
$\delta'_\mu\leq 0$ for all $\mu$.

Thanks to the estimate
\[
	\delta_\mu'-
	\inf_{\substack{v\in \mathcal{L}_\Gamma\\\mu(v)=1}} \rho(v)
	\leq
	\delta_\mu \leq
	\delta_\mu'+
	\sup_{\substack{v\in \mathcal{L}_\Gamma\\\mu(v)=1}} \rho(v)
\]
the modified and the usual critical exponents can be compared.
This leads to many applications involving the critical exponents
such as bounds on the limit set.
For example, one can apply \cite[Theorem~7.1]{DKO} to an $I$-Anosov subgroup $\Gamma$
with $I$ as in Corollary~\ref{cor:Anosovintro}.
One obtains
\[
	\dim \Lambda_I \leq \max_{\alpha\in I} \sup_{v\in \mathcal{L}_\Gamma} \frac{\rho(v)}{\alpha(v)},
\]
where $\Lambda_I$ is the limit set of $\Gamma$ in $G/P_I$ and
$P_I$ is the parabolic subgroup associated to $I$.

\subsubsection{Property (T)}
If $G$ has no factor locally isomorphic to $\mathfrak{so}(n,1)$ or $\mathfrak{su}(n,1)$
with $n\geq 2$,
then $G$ has Property (T).
This also leads to a bound on the growth indicator function.
\begin{proposition*}
	[{\cite[Theorem~7.1]{OhDichotomy}}]
	If $G$ is simple of higher rank and $\Gamma$ has infinite covolume,
	then
	\begin{equation}
		\label{eq:propTpsi}
		\psi_\Gamma \leq 2\rho -\Theta,
	\end{equation}
	where $\Theta$ is the half sum of all roots in a maximal strongly orthogonal system.
\end{proposition*}
We note that in general $\Theta < \rho$ as we sum only a proper subset of positive roots
	(see \cite{OhUniformPointwise} for precise values of $\Theta$).
Hence, the conclusion $\psi_\Gamma \leq \rho$
is much stronger than the $\Gamma$-independent bound $\psi_\Gamma \leq 2\rho - \Theta$.
If $\mathcal{L}_\Gamma$ avoids one facet of $\mathfrak{a}_+$,
then we can use \eqref{eq:propTpsi} together with Corollary~\ref{thm:onewallavoided}
to get an improved bound which is again independent of $\Gamma$.
This will be presented in Proposition~\ref{prop:boundwallavoided}.

\subsubsection{Examples}
\label{ss:nontempex}
Theorem~\ref{mainintro} provides a criterion when $\psi_\Gamma \leq \rho$ holds.
As noted above this is equivalent to the temperedness of $L^2(\Gamma \bk G)$.
Clearly, if $H$ is a reductive subgroup of $G$ such that 
$L^2(H \bk G)$ is non-tempered\footnote{This property is characterized by 
Benoist and Kobayashi \cite{BK1,BK2,BK3,BK4}
and extended by Benoist and Liang \cite{BenoistLiang}}
and $\Gamma$ is a lattice in $H$,
then $L^2(\Gamma \bk G)$ is non-tempered and $\psi_\Gamma\leq \rho$ does not hold.
By construction, $\Gamma$ in not Zariski dense in $G$
and it was unknown for some time whether there is a Zariski dense
discrete subgroup $\Gamma$ with $\psi_\Gamma \not \leq \rho$.
The first example has been discovered by Fraczyk and Oh \cite{FO25}
and is constructed by a deformation of the above procedure with $G=\mathrm{SO}_0(2,n)$ and $H=\mathrm{SO}_0(1,n)$.
For this example this paper gives strong bounds of $\psi_\Gamma$
and determines the maximal growth direction
under mild assumptions.
We refer to \cite[Section~6.3]{LWW} to see how the results are applied to
the example of \cite{FO25}.

The examples by Fraczyk and Oh are by construction $\{\alpha\}$-Anosov for a single simple root $\alpha$.
In contrast, there are also examples of $\{\alpha\}$-Anosov subgroups $\Gamma$
in higher rank
for which $\psi_\Gamma \leq \rho$ does hold
although Corollary~\ref{cor:Anosovintro} cannot be applied.
For example, if $G=\mathrm{PSL}_n(\R)$ and $\Gamma$ is a $(1,1,2)$-hyperconvex subgroup
whose Gromov boundary is a circle
or the image of a purely hyperbolic Schottky representation of a free group \cite{BT22},
then $\psi_\Gamma \leq \rho$.
These examples rely on calculation of the limit set for those examples
by Pozzetti-Sambarino-Wienhard \cite{PSW23} and Canary-Zhang-Zimmer \cite{CZZ22}
together with a bound of $\psi_\Gamma$ in terms of the dimension of the limit set \cite{DKO}.

\subsection{Idea of the proofs and outline of the article}
\label{sec:outline}
The proofs rely on the spectral results obtained in \cite{LWW}.
Let us describe the main ideas.
One can define a subset $\wt \sigma_\Gamma \subseteq \mathfrak{a}_\C^\ast$,
called the \emph{joint spectrum},
related to the spherical part of the decomposition of $L^2(\Gamma\bk G)$ into irreducible representations
(see Section~\ref{sec:spherical_dual}).
Using this, one defines for $\mu\in \mathfrak{a}_+^\ast$ the number
\[
	\theta_\mu \coloneqq \min \{t\geq 0\colon\Re  \wt \sigma_\Gamma \subseteq t\conv{\mu}\}.
\]
Since $\Re \wt \sigma_\Gamma \subseteq \theta_\mu \conv \mu$ for all $\mu\in \mathfrak{a}_+^\ast$
and the intersection of these convex hulls is again a convex hull of the same form,
we can define a critical element $\mu_\Gamma \in \mathfrak{a}_+^\ast$
with the property that
\[
	\theta_\mu  = \max_{v\in \mathfrak{a}_+} \frac{\mu_\Gamma(v)}{\mu(v)}
\]
for all $\mu\in \mathfrak{a}_+^\ast$.
We recall this construction and give alternative descriptions in Section~\ref{sec:defmugamma}.
In particular, we show that $\mu_\Gamma$ is given by the direction of maximal growth $v_\Gamma'$ in
Proposition~\ref{prop:mu_is_biggest_growth} if $\Gamma$ is Zariski dense
and prove Theorem~\ref{thm:spectralintro}.

One main result in \cite{LWW} is the following.
\begin{theorem*}
	[{\cite[Theorem~1.1]{LWW}}]
	\begin{equation}
	\label{thm:main}
		\forall \, \mu\in \mathfrak{a}_+^{\ast}\text{ with }\iota \mu = \mu
		\colon
		\quad
		\theta_\mu = \max(0,\delta'_\mu).
	\end{equation}
\end{theorem*}
If $\psi_\Gamma \not \leq \rho$,
\eqref{thm:main} implies that
\begin{equation*}
	\label{eq:thetabygif}
	\max_{v\in \mathcal{L}_\Gamma} \frac{(\psi_\Gamma-\rho)(v)}{\mu(v)} =\theta_\mu
\end{equation*}
for every $\mu\in \mathfrak{a}_+^\ast$
with $\iota\mu = \mu$ and consequently
(see Proposition~\ref{prop:exmax}),
\[
	\max_{v\in \mathfrak{a}_+} \frac{\mu_\Gamma(v)}{\mu(v)} = \theta_\mu
	=\max_{v\in \mathcal{L}_\Gamma} \frac{(\psi_\Gamma-\rho)(v)}{\mu(v)}
	\leq\max_{v\in \mathcal{L}_\Gamma} \frac{\mu_\Gamma(v)}{\mu(v)}.
\]
We conclude that
\[
	\max_{v\in \mathfrak{a}_+} \frac{\mu_\Gamma(v)}{\mu(v)}
\]
must be attained in $\mathcal{L}_\Gamma$.
	Depending on the assumption on $\mathcal{L}_\Gamma$ and $\mu_\Gamma$,
	we will apply this to different functionals $\mu$
	to obtain our results.
	Specifically, in Section~\ref{sec:fromlimitcone} we will make assumptions
	on the limit cone and derive properties of $\mu_\Gamma$ as well as bounds on $\psi_\Gamma$.
In particular, we prove Theorems~\ref{mainintro}, \ref{thm:oneweight}, and \ref{thm:psilinearintro}.

\subsection*{Acknowledgements}
I thank Tobias Weich for his help and advice in many aspects of this work.
I also thank Benjamin Delarue, Siwei Liang, Alex Nolte, and Hee Oh for discussions, helpful comments on the manuscript, and advice to the literature.
This work started within
the CRC “Integral Structures in Geometry and Representation Theory” Grant No. SFB-TRR 358/1 2023 - 491392403 
of the Deutsche Forschungsgemeinschaft (DFG)
and would not have been possible without its funding.

\section{Preliminaries}\label{sec:preliminaries}

\subsection{Notation}
\label{sec:not}
In this article, $G$ is a real semisimple connected non-compact Lie group with finite center
and $K$ is a maximal compact subgroup of $G$.
The corresponding Riemannian symmetric space of noncompact type is then $G/K$.
We fix an Iwasawa decomposition $G=KAN$, and have $A\cong \R^r$ where $r$ is the real rank of $G$ or the rank of the symmetric space $G/K$, respectively.
We denote by fractal letters
the corresponding Lie algebras.
For $g\in G$ let $H(g)\in \mathfrak{a}$ be the logarithm of the $A$-component in the Iwasawa decomposition.
Let $\Sigma \subseteq \mathfrak{a}^\ast$ be the root system of restricted roots,
$\Sigma^+$ the positive system corresponding to the Iwasawa decomposition,
$\Pi$ the subset of simple roots,
and $W$ the corresponding Weyl group acting on $\mathfrak{a}^\ast$.
As usual, for $\alpha\in \Sigma$, we denote
by $m_\alpha$ the dimension of the root space,
and by $\rho$ the half sum of positive restricted roots counted with multiplicity.
We denote by $(\omega_\alpha)_{\alpha\in \Pi} \subset \mathfrak{a}^\ast$ the (restricted) fundamental weights of $\Pi$
defined by
\[
	2 \frac{\langle \omega_\alpha,\beta\rangle}{\langle \beta,\beta\rangle}
	= \delta_{\alpha,\beta}, \quad \alpha,\beta\in \Pi,
\]
where  $\delta_{\alpha,\beta}$ is the Kronecker delta.
We note that the precise normalization of $\omega_\alpha$ does not matter
in the present article as we only consider the rays $\R_{\geq 0}
\omega_\alpha$.
Recall from the introduction that
$v_\alpha \in \mathfrak{a}_+  
=\{H\in \mathfrak{a}\mid \alpha(H)\geq 0 \:\forall \alpha\in \Sigma\}$
is defined by $\beta(v_\alpha)=0$ for $\beta\in \Pi\setminus \{\alpha\}$ and $\|v_\alpha\|=1$.
$v_\alpha$ and $\omega_\alpha$ correspond up to normalization.
Let $\mathfrak{a}_+^\ast$ the cone corresponding to $\mathfrak{a}_+$ in $\mathfrak{a}^\ast$
via the identification $\mathfrak{a} \leftrightarrow \mathfrak{a}^\ast$ through the Killing form $\langle\cdot,\cdot\rangle$.
Let us further introduce the opposition involution $\iota$:
There is a unique $w_0 \in W$ such that $w_0 (\Sigma^+) = -\Sigma^+$,
called the \emph{longest Weyl group element} and satisfies $w_0^2=1$.
The \emph{opposition involution} is defined as $\iota = -w_0 \colon \mathfrak{a}^\ast \to \mathfrak{a}^\ast$.
We define $\mathfrak{a}^{\ast,\mathrm{Her}} \coloneqq \ker (\iota -1)$
and $\mathfrak{a}^{\ast,\mathrm{Her}}_+ = \mathfrak{a}^{\ast,\mathrm{Her}}\cap \mathfrak{a}^\ast_+$.
For a non-empty subset $I$ of $\Pi$ let
$\mathfrak{a}_I\coloneqq \bigcap_{\alpha\in \Pi\setminus I} \ker \alpha$.

We fix a discrete subgroup $\Gamma\leq G$.
If we say $\Gamma$ is Zariski dense in $G$,
then we demand that $G$ is the identity component of the group of real points of a
semisimple algebraic group $\mathbf G$ defined over $\R$
and $\Gamma$ is dense with respect to the Zariski topology.
The only property that we need when we assume that $\Gamma$ is Zariski dense
is that the growth indicator function $\psi_\Gamma$ is concave.
Hence, everywhere in this article, the assumption \emph{$\Gamma$ is Zariski dense}
can be replaced by \emph{$\psi_\Gamma$ is concave}.

\subsection{The growth indicator function}
\label{ss:gif}

In this subsection we recall the definition of the growth indicator function $\psi_\Gamma$.
It was introduced by Quint \cite{Qui02}
generalizing the critical exponent in rank at least $2$.
For an open cone $\mathcal{C}\subseteq \mathfrak{a}$
let $\tau_{\mathcal{C}}$ be the abscissa of convergence for the series
$\sum_{\gamma\in \Gamma,\mu_+(\gamma)\in \mathcal{C}} e^{-s\|\mu_+(\gamma)\|}<\infty$,
	i.e. 
	\[
		\tau_{\mathcal{C}} = \inf\{s\in \R\mid 
		\sum_{\gamma\in \Gamma,\mu_+(\gamma)\in \mathcal{C}} e^{-s\|\mu_+(\gamma)\|}<\infty\}.
	\]
	The growth indicator function $\psi_\Gamma\colon \mathfrak{a}\to \R\cup \{-\infty\}$ is then defined as
	$\psi_\Gamma (u)= \|u\| \inf_{\mathcal{C}\ni u} \tau_{\mathcal{C}}$,
	where the infimum runs over all open cones $\mathcal{C}\subseteq \mathfrak{a}$
	containing $u$.
	We also set $\psi_\Gamma(0)=0$.
	One observes that $\psi_\Gamma$ is a positively homogeneous function
	that is upper semicontinuous.
	Moreover, $\psi_\Gamma$ is independent of the norm used on $\mathfrak{a}$.
	However, one usually uses the norm induced by the Killing form
	as it has the advantage of being invariant under the Weyl group
	which implies that $\psi_\Gamma$ is invariant under the opposition involution
	of $\mathfrak{a}$.
	One also finds that $\psi_\Gamma\leq 2\rho$,
	$\psi_\Gamma=-\infty$ outside $\mathfrak{a}_+$,
	and $\psi_\Gamma>-\infty$ implies $\psi_\Gamma \geq 0$.
	If $\Gamma$ is Zariski dense in a real algebraic group $G$,
	then one can make this more precise.
	Namely \cite{Qui02}, $\{v\in \mathfrak{a}\mid \psi_\Gamma(v)>-\infty\}$
	is precisely the limit cone
	\[
		\mathcal{L}_\Gamma = \{\lim_{i\to \infty} t_i \mu_+(\gamma_i) \in \mathfrak{a}_+ \mid
		t_i \to 0, \gamma_i\in \Gamma\}.
	\]
	Moreover, $\psi_\Gamma>0$ on the interior of $\mathcal{L}_\Gamma$
	which is non-empty \cite{Ben96}
	and $\psi_\Gamma$ is concave.

\subsection{Spherical dual, joint spectrum, and temperedness}\label{sec:spherical_dual}
\subsubsection{Spherical dual}
Let us denote with $\widehat G$ the unitary dual of $G$ and with $\widehat
G_{\mathrm{sph}}\subset \widehat G$ the spherical dual of $G$, i.e.~the set of
equivalence classes of irreducible unitary representations containing a
non-zero $K$-invariant vector.

$\widehat G_{\mathrm{sph}}$ can be parameterized by a $W$-invariant subset
of $\mathfrak{a}_\C^\ast$ or by picking a dominant element by a subset of $\mathfrak{a}_+^\ast + i \mathfrak{a}^\ast$ (see \cite[Thm.~IV.3.7]{gaga} for details).
The correspondence is as follows.
Let $\phi_\lambda$ be the elementary spherical function for $\lambda\in \mathfrak{a}_\C^\ast$,
i.e. \[
	\phi_\lambda = \int_{K}^{} {e^{-(\lambda+\rho)H(g^{-1}k)}} \: dk
	,
\]
where %
$H\colon G\to \mathfrak a$ is defined by $g\in Ke^{H(g)}N$.
For $\pi\in \widehat G_{\mathrm{sph}}$ let $v_K$ be a normalized $K$-invariant vector.
Then the function $\phi\colon G\to \C, \:\phi(g)=\langle \pi(g)v_K,v_K\rangle$
is an elementary spherical function $\phi_\lambda$ for some $\lambda\in \mathfrak{a}_\C$.
Moreover, $\phi_\lambda$ is positive definite as a diagonal matrix coefficient.
The induced map
$\widehat G_{\mathrm{sph}} \ni \mu\mapsto \lambda \in\{ \lambda\in \mathfrak{a}_\C^\ast \colon \phi_\lambda \text{ is positive definite}\}$ is a bijection.
The right hand side is $W$-invariant as $\phi_\lambda =\phi_\mu$ if and only if $\lambda\in W\mu$.

Every positive definite elementary spherical function $\phi_\lambda$ 
is Hermitian,
i.e.~it satisfies
$\phi_\lambda(g^{-1}) = \overline{\phi_\lambda(g)}$.
As $\phi_\lambda(g^{-1})=\phi_{-\lambda}(g)$ and $\overline{\phi_\lambda(g) }=\phi_{\overline \lambda}(g)$,
we must have $W(-\lambda)=W\overline{\lambda}$.
Hence, $\widehat G_{\mathrm{sph}} \subseteq \{\lambda\in \mathfrak{a}_\C^\ast \mid \exists w\in W \colon w\lambda =-\overline{ \lambda}\}$.
This implies $\Re \widehat G_{\mathrm{sph}}\cap \mathfrak{a}^\ast_+  \subseteq \mathfrak{a}^{\ast,\mathrm{Her}}$.

	We also recall that $\operatorname{conv}(W\mu)$ for $\mu\in \mathfrak{a}^\ast$
is the convex hull of the Weyl orbit $W\mu$ of $\mu$
which can be characterized by (see \cite[Lemma~IV.8.3]{gaga})
\begin{equation}
	\label{eq:convhull}
	\operatorname{conv}(W\mu)= \{\lambda\in \mathfrak{a}^\ast\mid \lambda(wH)\leq \mu(H) \:\forall H\in \mathfrak{a}_+, w\in W\}
\end{equation}
if $\mu\in \mathfrak{a}_+^\ast$.
Since every positive definite function on $G$ is bounded by its value at $1$,
$\widehat G_{\mathrm{sph}} \subseteq \operatorname{conv}(W\rho)+i\mathfrak{a}^\ast$ by \cite[Thm.~IV.8.1]{gaga}.

\subsubsection{Joint spectrum}
\label{ssec:jointspec}
Let us now give the definition of the joint spectrum for the discrete group $\Gamma \leq G$:
Consider the unitary representation $R$ on $L^2(\Gamma\backslash G)$ by right multiplication.
By the abstract Plancherel theory, it can be decomposed into a direct integral of irreducible representations
\[
	(R,L^2(\Gamma \backslash G))\simeq \int_{X}^{\oplus} {\pi_x} \: d\mu(x)
\]
where $(X,\mu)$ is a measure space and $\pi: X\longrightarrow \widehat G,\:
 x\longmapsto \pi_x$ 
 is a measurable map.

The \emph{joint spectrum} is now defined as
$\widetilde \sigma_\Gamma \coloneqq \supp (\pi_\ast \mu)\cap{\widehat G_{\mathrm{sph}}}
\subseteq \widehat G_{\mathrm{sph}}\subset \mathfrak{a}_\C^\ast.$

Let $\mathbb{D}(G/K)$ denote the algebra of $G$-invariant differential operators on $G/K$.
The \emph{Harish-Chandra isomorphism} (see \cite[Thm.~II.5.18]{gaga}) is the map
\[
	\chi: \left\{\begin{array}{rcl}
			\mathbb{D}(G/K) & \longrightarrow &\mathrm{Poly}(\mathfrak{a}_\C^\ast)^W\\
D& \longmapsto & \{\lambda\mapsto \chi_\lambda(D), \lambda\in \mathfrak a_\C^\ast\}
	             \end{array}
	\right.
\]
which is an algebra isomorphism between $\mathbb D(G/K)$ and the algebra of $W$-invariant polynomials on $\mathfrak a^*_\C$.
In particular, $\mathbb D(G/K)$ is abelian and is generated by $\mathrm{rank}(G/K)$ algebraically independent generators.
$\mathbb{D}(G/K)$ always contains the Laplace-Beltrami operator $\Delta$
and $\chi_\lambda(\Delta) = \|\rho\|^2 - \|\Re \lambda\|^2 + \|\Im \lambda\|^2$.
By $G$-invariance, each $D\in \mathbb{D}(G/K)$ acts as a differential operator
on $\Gamma\bk G/K$
and we can consider its spectrum $\sigma_{L^2(\Gamma\bk G/K)}(D)$.
We then have \[
 \widetilde \sigma_\Gamma = \{\lambda \in \mathfrak{a}_\C^\ast\mid \chi_\lambda(D) \in
	\sigma_{L^2(\Gamma\bk G/K)}(D) \quad \forall D\in \mathbb{D}(G/K)\}  \subset \mathfrak a^*_\C.
\]
and
\begin{equation}
	\label{eq:DGKspectrum}
	\sigma_{L^2(\Gamma\bk G/K)}(D)= \overline{\{\chi_\lambda(D) \colon \lambda\in \wt \sigma_\Gamma \}}
	\quad \text{for any}\quad  D\in \mathbb{D}(G/K) 
\end{equation}
by \cite[Section~3]{WW23}.
\subsubsection{Temperedness}\label{sec:tempered}
Temperedness is a well-studied property of unitary representations of $G$.
We recall its definition.
\begin{definition}
	[{\cite[Thm.~1 and 2]{CHH88}}]\label{prop:tempered_equivalences}
	A unitary $G$-representation
$(\rho, \mathcal  H)$ is called \emph{tempered}
if it has the following equivalent properties.
\begin{enumerate}
	\item $\pi$ is almost $L^2$,
		i.e.~there is a dense
subspace $V\subset\mathcal  H$  such that for any $v,w\in V$,
the matrix coefficient $g\mapsto \langle \rho(g) v,w\rangle$ lies in $L^q(G)$ for all $q>2$.
\item $\rho$ is weakly contained in the regular representation $L^2(G)$.
	\item For any $K$-finite unit vectors $v,w\in \mathcal  H$,
	\[
		|\langle \rho(g) v,w\rangle|\leq \left(\dim\langle Kv \rangle \dim \langle Kw\rangle\right)^{1/2}
		\phi_0(g)
	\]
	for any $g\in G$, where $\langle Kv\rangle$ denotes the subspace spanned by $Kv$.
\item If $(\rho,\mathcal{H})$ is decomposed as $\int_{X}^{\oplus} {\rho_x} \: d{\mu(x)} $ 
	then $\rho_x$ is weakly contained in $L^2(G)$ for $\mu$-almost every $x\in X$.
\end{enumerate}
\end{definition}
We call $\Gamma$ or $\Gamma\bk G$ tempered,
if the quasi-regular representation $L^2(\Gamma\bk G)$ is tempered.

\section{Definition of \texorpdfstring{$\mu_\Gamma$}{mu Gamma}}
\label{sec:defmugamma}
\subsection{Growth rates}
In this section we will give different definitions for the critical functional
$\mu_\Gamma$
and we will show that they are equivalent.
We will start with working directly with the discrete group $\Gamma$ and define
$\mu_\Gamma$ from the critical exponents.
Before doing so,
let us introduce the dual limit cone
\[
	\mathcal{L}_\Gamma^\star \coloneqq \{\mu\in \mathfrak{a}^\ast 
	\colon \mu(\mathcal{L}_\Gamma) \geq 0\}.
\]
let us recall the definition of the modified critical exponents
as stated in the introduction.
\[
\delta'_\mu = \inf\{s\in \R\colon \psi_\Gamma < s\mu +\rho\}
	= \sup_{v\in \mathcal{L}_\Gamma} \frac{(\psi_\Gamma-\rho)(v)}{\mu(v)}
\]
for $\mu\in \mathcal{L}_\Gamma^\star$.
Since in our setting $\psi_\Gamma$ often occurs with the $\rho$-shift,
we define
\[
	\psi_\Gamma' \coloneqq \psi_\Gamma - \rho .
\]
We also define the corresponding \emph{modified limit cone $\mathcal{L}_\Gamma'$} by
\[
	\mathcal{L}_\Gamma' = \{v \in \mathfrak{a} \colon
	\psi_\Gamma'(v) > 0\}
\]
which is an open subcone of $\mathcal{L}_\Gamma$.

We recall that $\delta'_\mu = \inf\{s\in \R \colon \sum_{\gamma\in \Gamma}	e^{-(s\mu+\rho)(\mu_+(\gamma))}<\infty\}$
if $\Gamma$ is Zariski dense
or if $\mu\in \operatorname{int} \mathcal{L}_\Gamma^\star = \{\mu\in \mathfrak{a}^\ast \colon \mu(\mathcal{L}_\Gamma \setminus \{0\}) >0\}$.

We use the following lemma.
\begin{lemma}
	If $\psi_\Gamma \not \leq \rho$, then the mapping
	\[
		\mathcal{L}_\Gamma^\star \setminus \{0\} \to \R\cup\{\infty\}, \mu \mapsto \|\mu\| \delta_\mu'
	\]
	is lower semicontinuous,
	positively homogeneous of order 0, and strictly convex where it is finite.
\end{lemma}
\begin{proof}
	By definition, $\delta_{t\mu}'= \frac 1t \delta_\mu'$ for $t>0$.
	This shows positive homogeneity of order $0$.

	Let $\mu\in \mathcal{L}_\Gamma^\star \setminus \{0\}$.
	By definition $\delta'_\mu = \inf\{t\geq 0 \colon \psi_\Gamma' \leq t\mu\}
	=\sup_{v\in \mathfrak{a}_+} \frac{\psi_\Gamma'(v)}{\mu(v)}$.
	By semicontinuity of $\psi_\Gamma'$,
	there exists $v\in \mathfrak{a}_+$
	such that $\psi_\Gamma'(v) > 0$ and
	$\frac{\psi_\Gamma'(v)}{\mu(v)}= \delta'_\mu $.
	For $\tilde \mu$ in a sufficiently small neighborhood of $\mu$,
$|\tilde \mu(v) -\mu(v)| \leq \varepsilon \|v\|$.
If $\mu(v)=0$, then $\delta'_{\tilde \mu} \geq \frac{\psi'_\Gamma(v)}{\tilde \mu(v)}
\geq \frac{\psi_\Gamma'(v)}{\varepsilon \|v\|}$.
Hence, $\delta'_{\tilde \mu}$ gets arbitrarily large when $\tilde \mu \to \mu$.
If $\mu(v) > 0$, then
$\tilde \mu(v) \geq  \mu(v) - \varepsilon \|v\|> 0$ for $\varepsilon$ sufficiently small.
Therefore,
$\delta'_{\tilde \mu} \geq
\frac{\psi'_\Gamma(v)}{\tilde \mu(v)} =\delta_\mu' \frac{\mu(v)}{\tilde \mu(v)}
\geq \delta'_\mu \left(1-\varepsilon \frac{\| v\|}{\mu(v) -\varepsilon \| v\|}\right)$.
Since the term in brackets gets close to $1$ as $\varepsilon \to 0$,
semicontinuity is proven.

For strict convexity, fix $\mu_1, \mu_2 \in \mathcal{L}_\Gamma^\star \setminus\{0\}$ linear independent
with $\delta'_{\mu_1}, \delta_{\mu_2}' < \infty$
as well as $0<s<1$.
We have $\psi_\Gamma' \leq \delta'_{\mu_i} \mu_i$ and clearly
$\psi_\Gamma' \leq t\delta'_{\mu_1} \mu_1 + (1-t) \delta'_{\mu_2}\mu_2$
for every $t \in [0,1]$.
Letting $t=\left(1+ \frac{(1-s)\delta_{\mu_1}'}{s\delta_{\mu_2}'}\right)^{-1}$
we have $t \frac{\delta'_{\mu_1}}{s} = (1-t)\frac{\delta'_{\mu_2}}{1-s} = \left(\frac{s}{\delta'_{\mu_1}} + \frac{1-s}{\delta'_{\mu_2}}\right)^{-1} \eqqcolon c$.
We thus find $\psi'_\Gamma \leq c (s\mu_1 +(1-s)\mu_2)$
and $\delta'_{s\mu_1 + (1-s)\mu_2} \leq c$.
To finish the proof we must show $c\|s\mu_1 + (1-s) \mu_2\| < \delta'_{\mu_1} \|\mu_1\| + \delta'_{\mu_2}\|\mu_2\|$.
Multiplying by $c$ and squaring gives
\[
	s^2\|\mu_1\|^2 + (1-s)^2 \|\mu_2\|^2 +2s(1-s) \langle\mu_1, \mu_2\rangle <
	\left(\left(s+(1-s)\frac{\delta'_{\mu_1}}{\delta'_{\mu_2}}\right) \|\mu_1\| +
	\left((1-s)+s\frac{\delta'_{\mu_2}}{\delta'_{\mu_1}}\right) \|\mu_2\|\right)^2.
\]
Since $\delta'_{\mu_i} >0$ and $s\in (0,1)$,
the right hand side is bigger than $(s\|\mu_1\| + (1-s)\|\mu_2\|)^2$.
Applying the Cauchy-Schwartz inequality $\langle\mu_1,\mu_2\rangle < \|\mu_1\| \|\mu_2\|$
finishes the proof.
\end{proof}
If $\psi_\Gamma \not \leq \rho$, by semicontinuity there is a ray $\R_{\geq 0} \mu$
where the infimum is attained.
By strict convexity, this ray is unique.
We can therefore define
\[
	\mu_\Gamma \coloneqq \delta'_\mu \mu
\]
which only depends on the ray $\R_{\geq 0} \mu$.
If $\psi_\Gamma \leq \rho$, we define $\mu_\Gamma=0$.
\begin{lemma}
	$\mu_\Gamma$ is invariant by the opposition $\iota$ and
	contained in the positive Weyl chamber,
	i.e.
	\[
		\mu_\Gamma \in \mathfrak{a}^{\ast,\mathrm{Her}}_+.
	\]
\end{lemma}
\begin{proof}
	If $\mu_\Gamma \notin \mathfrak{a}^\ast_+$,
	then there is $w\in W$ such that $w\mu_\Gamma = \mu'\in \mathfrak{a}^\ast_+$.
	But then $\mu_\Gamma(v) \leq \mu'(v)$ for all $v\in \mathfrak{a}_+$.
	This means $\psi'_\Gamma \leq \mu_\Gamma \leq \mu'$.
	Hence, $\delta'_{\mu'}\leq 1$
	and therefore $\delta_{\mu_\Gamma} \|\mu_\Gamma\| \geq \delta'_{\mu'} \|\mu'\|$.
	By strict convexity, this forces $ \mu_\Gamma =\mu' \in \mathfrak{a}_+^\ast$.

	To show that $\iota\mu_\Gamma =\mu_\Gamma$, we assume that they are different.
	As $\psi'_\Gamma$ is $\iota$-invariant,
	\[
		\psi'_\Gamma = t\psi'_\Gamma + (1-t) \psi'_\Gamma \circ \iota
		\leq t \mu_\Gamma + (1-t) \iota \mu_\Gamma.
	\]
	Thus, $\delta'_{t\mu_\Gamma + (1-t)\mu_\Gamma} \leq 1$.
	For the norm it is clear that $\|t\mu_\Gamma +(1-t)\mu_\Gamma\| < \|\mu_\Gamma\|$
	for all $t\in (0,1)$.
	This is a contradiction to the minimality of $\delta'_\mu \|\mu\|$.
\end{proof}

\subsection{Spectral side}
\label{ss:defspectral}
On the spectral side, $\mu_\Gamma$ was defined in \cite[Prop.~5.1]{LWW}.
Namely,
$\mu_\Gamma$ is the unique element of $\mathfrak{a}_+^\ast$ such that
\[
	\conv{\mu_\Gamma}
	= \bigcap_{\substack{\mu\in \mathfrak{a}_+^\ast\\\Re \wt \sigma_\Gamma\subseteq \conv{\mu}}} 
	\conv{\mu}.
\]
Since for all $\mu\in \mathfrak{a}_+^\ast$,
$\Re \wt \sigma_\Gamma \subseteq \conv {\theta_\mu \mu}$ if $\theta_\mu$ defined in \eqref{thm:main} is finite,
we have $\conv{\mu_\Gamma} \subseteq \conv {\theta_\mu \mu}$.
By \cite[Lemma~IV.8.3]{gaga}, this relation of the convex hulls means
\begin{equation}
	\label{eq:unglconvexhulls}
	\mu_\Gamma(v)\leq \theta_\mu \mu(v)
\end{equation}
for all $v\in \mathfrak{a}_+$.
By minimality of $\theta_\mu$, we therefore have $\theta_\mu = \sup_{v\in \mathfrak{a}_+} \frac{\mu_\Gamma(v)}{\mu(v)}$.
Let us assume $\psi_\Gamma \not \leq \rho$.
Then $\theta_{\mu_\Gamma}=1$.
If $\mu_\Gamma$ and $\mu$ are linearly independent,
\[
	\|\mu_\Gamma\| \theta_{\mu_\Gamma} = \frac{\|\mu_\Gamma\|^2 \|\mu\|}{\|\mu_\Gamma\| \|\mu\|} < 
	\frac{\|\mu_\Gamma\|^2 \|\mu\|}{\langle \mu_\Gamma, \mu\rangle}
	\leq \sup_{v\in \mathfrak{a}_+}	\frac{\mu_\Gamma(v) \|\mu\|}{\mu(v)}
	=\|\mu\|\theta_\mu.
\]
Hence, $\R_{\geq 0} \mu_\Gamma$ minimizes $\|\mu\| \theta_\mu$.
By \eqref{thm:main}, $\theta_\mu=\delta_\mu'$ if $\psi_\Gamma \not \leq \rho$,
so that the two definitions coincide.
If $\psi_\Gamma \leq \rho$,
we have $\Re \wt \sigma_\Gamma \subseteq i \mathfrak{a}^\ast$
by \cite[Corollary~1.3]{LWW}.
Hence, $\mu_\Gamma=0$ as before.
\begin{remark}
	One might think that $\mu_\Gamma$ is uniquely determined
	by $\theta_\mu$ for $\mu=\omega_\alpha$, $\alpha\in \Pi$,
	being the extremal points of the cone $\mathfrak{a}_+^\ast$.
	This is however not the case a priori.
	To see this consider the following example.
	Let $\Sigma$ be of type $B_3$ which is e.g.~the case for $G=\mathrm{SO}(3,n)$, $n\geq 4$.
	Here, $\mathfrak{a}\simeq \mathfrak{a}^\ast\simeq \R^3$
	and
	\[
		\mathfrak{a}_+ =\{v=(v_1,v_2,v_3)\colon v_1\geq v_2 \geq v_3 \geq 0\}.
	\]
	Furthermore, we have
	\[
		\omega_1=e_1, \quad \omega_2=e_1+e_2,
		\quad\text{and}\quad \omega_3=\frac 12 (e_1+e_2+e_3)
	\]
	according to the simple roots
	\[
		\Pi= \{\alpha_1=e_1-e_2, \alpha_2= e_2-e_3,\alpha_3 = e_3\}.
		\]
	Let $\mu_\Gamma = \mu_1 e_1 + \mu_2 e_2 + \mu_3 e_3$ with $\mu_1\geq \mu_2 \geq \mu_3 \geq 0$.
	Then
	\[
		\theta_{\omega_1} = \max_{v\in \mathfrak{a}_+} \frac{\mu_\Gamma(v)}{\omega_1(v)}
		= \max(\mu_1,\mu_1+\mu_2,\mu_1+\mu_2+\mu_3) = \mu_1+\mu_2+\mu_3
	\]
	by considering the three rays $\R_{\geq 0} \omega_i$.
	Also
	\[
		\theta_{\omega_3} = \max\left(\mu_1, \frac 1{2}(\mu_1+\mu_2), \frac 13 (\mu_1+\mu_2+\mu_3)\right) = \mu_1
	\]
	and
	\[
		\theta_{\omega_2} = \max\left(\mu_1, \frac 1{2}(\mu_1+\mu_2), \frac 12 (\mu_1+\mu_2+\mu_3)\right) = 
		\begin{cases}
			\mu_1 &\colon \mu_1 \geq \mu_2 +\mu_3\\
			\frac 12 (\mu_1+\mu_2+\mu_3)&\colon \mu_1\leq \mu_2 +\mu_3
		\end{cases}
		.
	\]
Hence, we can only determine the sum $\mu_2+\mu_3$ from $\theta_{\omega_i}$
but not the two parameters individually.
\end{remark}

\subsection{Tent property}
In \cite[Thm.~2.5]{KMO24} it has been shown that
\begin{equation}
	\label{eq:tent}
	\psi_\Gamma \leq \inf_{\substack{\mu\in \mathcal{L}^\star_\Gamma
	\\ \delta_\mu <\infty}} \delta_\mu \mu
\end{equation}
under the assumption that $\Gamma$ is Zariski dense.
We recall that here
\[
	\delta_\mu \coloneqq\inf\left\{s\in \R \colon \sum_{\gamma\in \Gamma}
	e^{-s\mu(\mu_+(\gamma))}<\infty\right\}.
\]
The assumption of Zariski density
is used to push the estimate $\psi_\Gamma \leq \delta_\mu \mu$ from the interior of the limit cone
to its closure and therefore to all of $\mathfrak{a}_+$.
We note that their proof works the same way
without the assumption that $\Gamma$ is Zariski dense
as long as $\mu\in \operatorname{int}\mathcal{L}^\star_\Gamma$.
By our definition
\[
	\delta'_\mu 
	= \inf\{t\in \R\colon \psi'_\Gamma < t \mu\}
	= \sup_{v\in \mathfrak{a}_+} \frac{\psi'_\Gamma(v)}{\mu(v)}.
\]
Hence, we have $\psi'_\Gamma \leq \delta'_\mu \mu$ for all $\mu\in \mathcal{L}^\star_\Gamma$
with $\delta'_\mu <\infty$ by definition.
In particular, we also have
\begin{equation}
	\label{eq:tent'}
	\psi'_\Gamma \leq \inf_{\substack{\mu\in \mathcal{L}_\Gamma^\star
	\\ \delta'_\mu <\infty}} \delta'_\mu \mu.
\end{equation}
\cite[Theorem~1.1]{LWW} allows us to determine this infimum
when $\mu$ is restricted to the positive Weyl chamber.
\begin{proposition}
	\label{prop:tentattained}
	We have
	\begin{equation*}
	\psi'_\Gamma \leq \inf_{\substack{\mu\in \mathcal{L}_\Gamma^\star 
	\\ \delta'_\mu <\infty}} \delta'_\mu \mu
\end{equation*}
and,
if $\psi_\Gamma \not\leq \rho$,
\begin{equation*}
	\mu_\Gamma = \inf_{\substack{\mu\in \mathfrak{a}_+^\ast
	\\ \delta'_\mu <\infty}} \delta'_\mu \mu
	\quad \text{on}\quad \mathfrak{a}_+^{\mathrm{Her}}.
\end{equation*}
\end{proposition}
\begin{proof}
	By \eqref{eq:unglconvexhulls}, $\mu_\Gamma(v) \leq \theta_\mu \mu(v)$
	for all $\mu\in \mathfrak{a}^\ast_+$ and $v\in \mathfrak{a}_+$.
	By \eqref{thm:main}, $\theta_\mu = \delta'_\mu$ for all $\mu\in \mathfrak{a}_+^{\ast,\mathrm{Her}}$.
	It follows that
	\[
		\mu_\Gamma(v) \leq \inf_{\substack{\mu\in \mathfrak{a}_+^{\ast,\mathrm{Her}}
		\\ \delta'_\mu <\infty}} \delta'_\mu \mu(v)
	\]
	on $\mathfrak{a}_+$.
	Using the $\iota$-invariance of $\psi_\Gamma$,
	one easily finds that
	\[
		\delta'_{\frac 12(\mu+\iota\mu)} \leq \delta'_\mu = \delta'_{\iota\mu}.
	\]
	This implies
	\[
		\inf_{\substack{\mu\in \mathfrak{a}_+^{\ast,\mathrm{Her}}
		\\ \delta'_\mu <\infty}} \delta'_\mu \mu
		\leq \inf_{\substack{\mu\in \mathfrak{a}_+^{\ast}
		\\ \delta'_\mu <\infty}} \delta'_\mu \mu
	\]
	on $\mathfrak{a}_+^{\ast,\mathrm{Her}}$
	as
	\[
		\frac 12 (\mu+\iota\mu)(v) =
		\mu(v)
	\]
	for $v\in \mathfrak{a}_+^{\ast,\mathrm{Her}}$.
	This completes the proof.
\end{proof}

\subsection{Relation to direction of maximal growth}
Let us now relate $\mu_\Gamma$ with the direction of fastest growth.
To that end, let
\[
	\delta'\coloneqq \max_{\|v\|=1} \psi'_\Gamma(v)
\]
and choose $v_\Gamma\in \mathfrak{a}_+$
normalized with
\[
	\max(0,\delta') = \psi'_\Gamma(v_\Gamma).
\]
As usual we identify $\mathfrak{a}$ and $\mathfrak{a}^\ast$ by the given inner product
so that we can consider $\delta'_{v_\Gamma}$.
We observe that
\[
	\delta' = \psi'_\Gamma(v_\Gamma) \leq \sup_{\langle v, v_\Gamma\rangle =1} \psi'_\Gamma(v)=\delta'_{v_\Gamma}.
\]
Moreover, for any $\mu\in \mathfrak{a}_+^\ast \setminus \{0\}$ not proportional to $v_\Gamma$,
\begin{equation}
	\label{eq:exponents_growth}
	\|\mu\| \delta'_\mu = \|\mu\| \sup_{v} \frac{\psi'_\Gamma(v)}{\mu(v)} \geq
	\|\mu\| \frac{\psi'_\Gamma(v_\Gamma)}{\mu(v_\Gamma)} > \psi'_\Gamma(v_\Gamma)
	=\delta'.
\end{equation}
The following proposition shows that we have equality for $\mu = v_\Gamma$.
\begin{proposition}
	\label{prop:mu_is_biggest_growth}
	If $\Gamma$ is Zariski dense and $\psi_\Gamma \not \leq \rho$, then $v_\Gamma$ is unique,
	\[
		\delta' = \delta_{v_\Gamma}',
	\]
	and
	\[
		\mu_\Gamma = \delta'_{v_\Gamma} v_\Gamma.
	\]
\end{proposition}
\begin{proof}
	The uniqueness of $v_\Gamma$ directly follows from the concavity of $\psi'_\Gamma$
	and the strict convexity of the unit ball.
	If we can show $\delta'= \delta'_{v_\Gamma}$,
	then $\mu_\Gamma = \delta'_{v_\Gamma} v_\Gamma$ follows from
	$\|\mu_\Gamma\| \delta'_{\mu_\Gamma} = \min_\mu \|\mu\| \delta'_\mu$
	and \eqref{eq:exponents_growth}.

	To show this, we pick $v_n\in \mathfrak{a}_+$ with
	$\delta'_{v_\Gamma}=\lim_n \psi_\Gamma'(v_n)$
	and $1={\langle v_n, v_\Gamma \rangle}$.
	By concavity 
	$\psi'(tv_n + (1-t) v_\Gamma) \geq \psi'_\Gamma(v_\Gamma)
	+ t(\psi'_\Gamma(v_n) - \psi'_\Gamma(v_\Gamma))$
	for all $t\in [0,1]$.
	Letting $v_n(t)\coloneqq tv_n +(1-t)v_\Gamma$,
	we thus have $\frac 1t (\psi'_\Gamma(v_n(t)) - \psi'_\Gamma(v_\Gamma)) \geq
	\psi'_\Gamma(v_n) - \psi'_\Gamma(v_\Gamma)$.
	On the other hand,
	\[
		\psi'_\Gamma(v_n(t))-\psi'_\Gamma(v_\Gamma)
		=\|v_n(t)\| \psi'_\Gamma\left(\frac{v_n(t)}{\|v_n(t)\|}\right) - \psi'_\Gamma(v_\Gamma)
		\leq (\|v_n(t)\| -1) \psi'_\Gamma(v_\Gamma)
	\]
	by definition of $v_\Gamma$.
	Since
	\[
		\|v_n(t)\|^2 = t^2 \|v_n\| + 2t(1-t)\langle v_n, v_\Gamma \rangle + (1-t)^2 \|v_\Gamma\|^2
		= t^2 (\|v_n\|^2 -1) +1,
	\]
	we have
	\begin{align*}
		\lim_{t\to 0} \frac{\|v_n(t)\| -1}{t}
		&= \lim_{t\to 0} \frac{\sqrt{t^2(\|v_n\|^2 -1) +1} -1}{t}
		= \lim_{t\to 0} \frac{2t(\|v_n\|^2 -1)}{2\sqrt{t^2(\|v_n\|^2 -1) +1}}
		=0
	\end{align*}
	by l'Hôpital's rule.
	Therefore, $\psi'_\Gamma(v_\Gamma) \geq \psi'_\Gamma(v_n)$ for all $n$.
	This implies $\delta' \geq \delta'_{v_\Gamma}$.
\end{proof}

\begin{remark}
	If $\Gamma$ is $\Pi$-Anosov, then $\psi_\Gamma$ is strictly concave
	and differentiable
	and $v_\Gamma \in \operatorname{int}(\mathfrak{a}_+)$ \cite{Sam}.
In the proof of the proposition,
	we can choose $v_n=v_\infty$ and
	we get $\psi'_\Gamma(v_\Gamma) > \psi'_\Gamma(v_\infty)$.
	However, by definition of $\delta'_{v_\Gamma}$, $\psi'_\Gamma(v_\Gamma) \leq \psi'_\Gamma(v_\infty)$.
	It follows that $v_\Gamma = v_\infty$,
	i.e.~$\delta'_{v_\Gamma}$ is realized at $v_\Gamma$.
	This is part of \cite[Lemma~2.24]{ELO23}.
\end{remark}

Proposition~\ref{prop:mu_is_biggest_growth} allows us to apply the argument of \cite[Section~6.2]{LWW} in general.
This is the content of the following proposition.
	 \begin{proposition}
		 \label{prop:muGammainspectrum}
		If $\Gamma$ is Zariski dense,
		then
		\[
			\mu_\Gamma\in \wt \sigma_\Gamma
			\quad
			\text{and}
			\quad
			\chi_{\mu_\Gamma}(\Delta) = \min \sigma_{L^2(\Gamma\bk G/K)}(\Delta).
		\]
	 \end{proposition}
	 \begin{proof}
		By \eqref{eq:DGKspectrum}, there exist
	 	let $\mu_n\in \wt \sigma_\Gamma$ such that
		$\chi_{\mu_n}(\Delta)\to\min \sigma(\Delta)$.
		By \cite{WZ23},
		\[
		\min \sigma(\Delta) = 
		\| \rho\|^2 -  \max\left(0, \sup_{v\in \mathfrak{a}_+}
	\frac{\psi_\Gamma - \rho(v)}{\|v\|}\right)^2
= \|\rho\|^2 - \max (0, \delta')^2
\]
	Therefore, 
	as
	$\chi_{\mu_n}(\Delta) =\|\rho\|^2 -\|\Re \mu_n\|^2 + \|\Im \mu_n\|^2$,
	\[
		\|\Re \mu_n\|^2 - \|\Im \mu_n\|^2\to \max(0, \delta')^2.
\]
We have $\|\Re \mu_n\|^2 \leq \|\mu_\Gamma\|^2$, as $\operatorname{conv}(W\mu_\Gamma)$ is contained in the ball of radius $\|\mu_\Gamma\|$.
If $\psi_\Gamma \not \leq \rho$, $0< \delta' = \delta'_{v_\Gamma} = \|\mu_\Gamma\|^2$
by Proposition~\ref{prop:mu_is_biggest_growth}.
Hence,
\[
	\|\mu_\Gamma\|^2 -\|\Im \mu_n\|^2 \geq \|\Re \mu_n\|^2 - \|\Im \mu_n\|^2 \to \|\mu_\Gamma\|^2
\]
and therefore $\Im \mu_n \to 0$ as well as $\|\Re \mu_n \|\to \|\mu_\Gamma\|$.
We infer that $\Re \mu_n \to \mu_\Gamma$ as
\[
	\conv{\mu_\Gamma} \cap \{\lambda \in \mathfrak{a}^\ast\mid \|\lambda\|=\|\mu\|\}=W\mu_\Gamma
\]
by strict convexity of the sphere.
To conclude, we use that $\wt \sigma_\Gamma$ is closed and $\mu_n \to \mu_\Gamma$.
If $\psi_\Gamma \leq \rho$, $\Re \mu_n=0$ and $\delta'\leq 0$.
This implies $\Im \mu_n=0$, so that $\mu_n=0=\mu_\Gamma$.
	 \end{proof}
	 
	 From this we can determine the function $\sup_{\lambda\in \wt \sigma_\Gamma} \Re\lambda$ used intensively in \cite{LWW}.
	 \begin{corollary}
	 	Assume that $\Gamma$ is Zariski dense.
		Then
		\[
			\sup_{\lambda\in \wt \sigma_\Gamma} \Re \lambda = \mu_\Gamma
			\quad
		\text{on}
		\quad
			\mathfrak{a}_+.
		\]
	 \end{corollary}
	 
	 \begin{proof}
	 	By definition of $\mu_\Gamma$,
		$\Re \wt \sigma_\Gamma \subseteq \conv{\mu_\Gamma}$
		which means $\Re\lambda (v) \leq \mu_\Gamma(v)$ for all $v\in \mathfrak{a}_+$.
		Hence, $\sup_{\lambda\in \wt \sigma_\Gamma} \Re \lambda(v) \leq \mu_\Gamma(v)$.
		We must have equality as $\mu_\Gamma \in \wt \sigma_\Gamma$.
	 \end{proof}

	 \begin{proof}
		 [Proof of Theorem~\ref{thm:spectralintro}]
		 Combine this corollary and Proposition~\ref{prop:muGammainspectrum}
		 with Proposition~\ref{prop:mu_is_biggest_growth}.
	 \end{proof}

\section{Proofs of the main results}
\label{sec:fromlimitcone}
In this section we obtain results under assumptions on the position of the limit cone
in $\mathfrak{a}_+$.
In particular, we will prove Theorem~\ref{mainintro}.

	The proofs follow similar arguments as in \cite[Thm.~1.4]{LWW}
	but choosing a specific functional.
	Let us therefore summarize these lines of argument.
	\begin{proposition}
		\label{prop:exmax}
		Assume $\psi_\Gamma \not \leq \rho$.
		Then for any $\mu
 		\in  \mathfrak{a}^{\ast,\mathrm{Her}}_+$,
		there is $v_0 \in \mathcal{L}_\Gamma'$ such that
		\[
			\sup_{v\in \mathfrak{a}_+} \frac{\mu_\Gamma(v)}{\mu(v)} = \frac{\psi'_\Gamma(v_0)}{\mu(v_0)}.
		\]
		In particular,
		\[
			\sup_{v\in \mathfrak{a}_+} \frac{\mu_\Gamma(v)}{\mu(v)} = \frac{\mu_\Gamma(v_0)}{\mu(v_0)}.
		\]
	\end{proposition}
	We note that this equality also holds if $\delta'_{\mu}$ is infinite.
	In this case $\mu(v_0)=0$
	with the convention $0/0 =0$.
	\begin{proof}
		If $\mu_\Gamma=0$ there is nothing to show.
    Let us assume that $\mu_\Gamma\neq 0$.
    We observed in Section~\ref{sec:defmugamma}
    that $\theta_\mu = \sup_{v\in \mathfrak{a}_+} \frac{\mu_\Gamma(v)}{\mu(v)}$
    which is $>0$ as $\mu_\Gamma \neq 0$.
    By \eqref{thm:main},
    \[
	    0<\theta_\mu = \delta'_\mu =
	    \inf \{t\in \R \colon \psi'_\Gamma < t\mu\}
	    =\sup_{v\in \mathfrak{a}_+} \frac{\psi'_\Gamma(v)}{\mu(v)}.
    \]
    Since $\psi'_\Gamma$ is upper semicontinuous,
    there is $v_0\in \mathfrak{a}_+$ such that $\theta_\mu = \frac{\psi'_\Gamma(v_0)}{\mu(v_0)}.$
    Clearly, $v_0 \in \mathcal{L}_\Gamma'$ by definition of $\mathcal{L}'_\Gamma$.
    For the 'in particular' part, it suffices to observe that $\psi'_\Gamma \leq \mu_\Gamma$.
	\end{proof}
The next theorem is the main result of this section.	
	
 	\begin{theorem}
		\label{thm:muregularimplieslimitcone}
 		If $\langle \mu_\Gamma, \alpha \rangle >0$,
		then there is non-zero $v\in C_\alpha=\R_{\geq 0} v_\alpha + \R_{\geq 0}v_{\iota\alpha}$
		such that
		\begin{equation}
			\label{eq:ofmainthm}
			\psi_\Gamma(v) = \mu_\Gamma(v) + \rho(v).
		\end{equation}
		In particular, $v\in \mathcal{L}'_\Gamma$.
 		If $\iota\alpha = \alpha$, then $v_\alpha \in \mathcal{L}'_\Gamma$
		and $\psi_\Gamma(v_\alpha)= \mu_\Gamma(v_\alpha) +\rho(v_\alpha)$.
 	\end{theorem}

 	\begin{proof}
 		For $t>0$ let $\lambda_t =\mu_\Gamma - t (\alpha+\iota \alpha)$.
 		Then $\langle \lambda_t, \alpha \rangle \geq 0$
		for $t\leq \frac{\langle\mu_\Gamma,\alpha\rangle}{\langle \alpha,\alpha + \iota\alpha\rangle}$.
 		As $\iota \mu_\Gamma =\mu_\Gamma$, $\langle\iota \alpha,\mu_\Gamma\rangle=\langle\alpha,\mu_\Gamma\rangle > 0$
 		and therefore $\langle \iota \alpha, \lambda_t\rangle \geq 0$
 		for such $t$.
 		Moreover, for $\beta\in \Pi \setminus \{\alpha,\iota\alpha\}$,
 		$\langle \lambda_t, \beta \rangle =
 		\langle \mu_\Gamma , \beta \rangle -t\langle \alpha,\beta\rangle
 		\geq \langle \mu_\Gamma , \beta \rangle\geq 0$
 		as $\langle \alpha,\beta\rangle \leq 0$
 		and $\langle \iota \alpha,\beta\rangle \leq 0$.
 		Hence, $\lambda_t\in \mathfrak{a}_+^\ast$.
 		The linear optimization problem $\max_{v\in \mathfrak{a}_+,\lambda_t(v)=1} \mu_\Gamma(v)$
 		is attained in one of the extremal points of the admissible region,
 		i.e.
 		\[
 			\max_{v\in \mathfrak{a}_+} \frac{\mu_\Gamma(v)}{\lambda_t(v)}
 				=\frac{ \mu_\Gamma(v_\beta)}
				{\lambda_t(v_\beta)}
 		\]
 		for some $\beta\in \Pi$.
 		By definition of $\lambda_t$, for $\beta\in \Pi\setminus\{\alpha,\iota \alpha\}$,
 		\[
			\frac{\mu_\Gamma( v_\beta)}{ \lambda_t(v_\beta)}
			=\frac{ \mu_\Gamma(v_\beta)}{ \mu_\Gamma(v_\beta)} =1
 		\]
 		and
 		\[
			\frac{\mu_\Gamma (v_{\iota\alpha})}{\lambda_t(v_{\iota\alpha}) }
 				=
			\frac{\mu_\Gamma (v_{\alpha})}{\lambda_t(v_{\alpha}) }
 				=
				\frac{ \mu_\Gamma(v_\alpha)}
				{ \mu_\Gamma(v_\alpha) -t (\alpha + \iota\alpha) (v_\alpha) }
 				>1
 			\]
			for $0<t< \frac{\mu_\Gamma(v_\alpha)}{(\alpha+\iota\alpha)(v_\alpha)}$.
			Lemma~\ref{la:posofweight} below
			shows $\mu_\Gamma(v_\alpha)>0$,
			so that we can fix such $t$.
			We can indeed take $t= \frac{\langle\mu_\Gamma,\alpha\rangle}{\langle \alpha,\alpha + \iota\alpha\rangle}$
			but this does not give a stronger result.

 			We note that
 			\begin{equation}
				\label{eq:setmaximum}
 				\left\{v\in \mathfrak{a}_+\colon
					\frac{\mu_\Gamma(v_\alpha)}
					{ \lambda_t(v_\alpha)}=
 					\frac{\mu_\Gamma(v)}{\lambda_t(v)}
 				\right\}
				=C_\alpha.
 			\end{equation}
			We now apply Proposition~\ref{prop:exmax}
			to obtain $v_0 \in \mathcal{L}'_\Gamma$
			such that
			\[
					\frac{\mu_\Gamma(v_\alpha)}
					{ \lambda_t(v_\alpha)}=
				\sup_{v\in \mathfrak{a}_+} \frac{\mu_\Gamma(v)}{\lambda_t(v)}
				= \frac{\psi_\Gamma'(v_0)}{\lambda_t(v_0)} = \frac{\mu_\Gamma(v_0)}{\lambda_t(v_0)}.
			\]
			We note that $v_0$ must be non vanishing.
			By \eqref{eq:setmaximum},
			we get
			$v_0\in C_\alpha$.
			It remains to observe that $\mu_\Gamma(v_0) >0$.
			This is clear as $\mu_\Gamma(v_\alpha)>0$ by Lemma~\ref{la:posofweight}.
 	\end{proof}

	\begin{lemma}
		\label{la:posofweight}
		For $\mu\in \mathfrak{a}^{\ast,\mathrm{Her}}_+$ and $\alpha\in \Pi$
		we have
		\[
			\langle \mu, \alpha \rangle
			\leq 	\frac{\mu(v_\alpha) }{ (\alpha+\iota\alpha)(v_\alpha)} 
			\langle \alpha, \alpha +\iota\alpha\rangle.
		\]

	\end{lemma}
	\begin{proof}
		We write $\mu = \sum_{\beta\in \Pi} c_\beta \beta$ with $c_\beta\in \R$.
		$\mu\in \mathfrak{a}^{\ast}_+$ implies $c_\beta \geq 0$.
		Then $\mu(v_\alpha) = \alpha(v_\alpha) c_\alpha$.
		By $\iota$-invariance of $\mu$ we also have
		\[
			\mu(v_\alpha)=\mu(v_{\iota\alpha})= \iota\alpha(v_{\iota\alpha}) c_ {\iota\alpha}=\alpha(v_\alpha)c_{\iota\alpha}.
		\]
		Thus $c_\alpha=c_{\iota\alpha}$.
		We now calculate
		for $\iota \alpha \neq \alpha$:
		\[
			\langle \mu, \alpha \rangle = c_\alpha \langle \alpha,\alpha\rangle
			+ c_{\iota\alpha} \langle \alpha, \iota\alpha\rangle
			+ \sum_{\alpha \neq \beta\neq \iota\alpha} c_\beta \langle \alpha, \beta\rangle
			\leq 
			c_\alpha \langle \alpha, \alpha +\iota\alpha\rangle
		\]
		as $\langle \alpha,\beta\rangle \leq 0$.
		We conclude $	\langle \mu, \alpha \rangle
		\leq \frac{\mu(v_\alpha) }{ \alpha(v_\alpha)} 
\langle \alpha, \alpha +\iota\alpha\rangle =
\frac{\mu(v_\alpha)}{ (\alpha + \iota\alpha)(v_\alpha)} 
\langle \alpha, \alpha +\iota\alpha\rangle.$

If $\iota\alpha = \alpha$,
then by the same arguments,
\[
	\langle \mu, \alpha \rangle \leq c_\alpha \langle \alpha,\alpha\rangle
	= \frac 12 c_\alpha \langle \alpha , \alpha+ \iota\alpha\rangle
	=\frac{\mu(v_\alpha)}{2\alpha(v_\alpha)}  \langle \alpha , \alpha+ \iota\alpha\rangle
	=\frac{\mu(v_\alpha)}{(\alpha + \iota\alpha)(v_\alpha)}  \langle \alpha , \alpha+ \iota\alpha\rangle.
\]
This completes the proof.
\end{proof}
\begin{remark}
	We note that the assumption $\langle \mu_\Gamma, \alpha\rangle >0$
	cannot be weakened as we describe below.
	Let us for simplicity assume that $\iota\alpha =\alpha$.
	In the proof of Theorem~\ref{thm:muregularimplieslimitcone}
	we used that there is $\lambda \in \mathfrak{a}_+^{\ast,\mathrm{Her}}$
	such that the maximum $\max_{v\in \mathfrak{a}_+} \frac{\mu_\Gamma(v)}{\lambda(v)}$
	is attained at a unique ray given by $\R_{\geq 0} v_\alpha$.
	The existence of such $\lambda$ is in fact equivalent to $\langle \mu_\Gamma,\alpha\rangle >0$.
	Indeed, for any $t \in [0,\frac 12]$,
	$\omega_\alpha -t\alpha \in \mathfrak{a}^\ast_+$
	as $\langle \omega_\alpha - t\alpha , \alpha\rangle = (\frac 12-t) \|\alpha\|^2 \geq 0$
	and $\langle \omega_\alpha - t\alpha, \beta\rangle
	=-t \langle\alpha,\beta \rangle \geq 0$.
	Hence,
	if $\langle \mu_\Gamma,\alpha\rangle$ were $0$,
	\[
		\frac{\langle\mu_\Gamma,\omega_\alpha - t \alpha\rangle}{\langle\lambda,\omega_\alpha - t\alpha\rangle}
		=\frac{\langle\mu_\Gamma,\omega_\alpha \rangle}{\langle\lambda,\omega_\alpha\rangle - t\langle\lambda,\alpha\rangle}
		\geq 
		\frac{\langle\mu_\Gamma,\omega_\alpha \rangle}{\langle\lambda,\omega_\alpha \rangle}
	\]
	since $\langle \lambda,\alpha\rangle \geq 0$.
	Thus, $\R_{\geq 0} v_\alpha$ (which corresponds to $\R_{\geq 0} \omega_\alpha$) is not the unique maximizer of $\max_{v\in \mathfrak{a}_+} \frac{\mu_\Gamma(v)}{\lambda(v)}$.
\end{remark}
The following corollary proves Theorem~\ref{mainintro}.
 	\begin{corollary}
		\label{cor:twowallsavoided}
		The following statements are equivalent:
		\begin{enumerate}
			\item $C_\alpha \cap \mathcal{L}_\Gamma' = \emptyset $
 		for all $\alpha \in \Pi$.
	\item $\mathcal{L}'_\Gamma =\emptyset$.
	\item $\mu_\Gamma = 0$.
	\item $\psi_\Gamma \leq \rho$.
	\item $\widetilde \sigma_\Gamma \subseteq i \mathfrak{a}^\ast$.
	\item $L^2(\Gamma \bk G)$ is tempered.
	\item $\min \sigma(\Delta)=\|\rho\|^2$.
	\item For all $\varepsilon >0$, there is $d_\varepsilon$ such that
		for all $f_1,f_2 \in L^2(\Gamma\bk G)^K$ and $v\in \mathfrak{a}_+$:
		\[
			|\langle (\exp v) \cdot f_1,f_2 \rangle_{L^2(\Gamma\bk G)}|
			\leq d_\varepsilon e^{\varepsilon \|v\|} e ^{-\rho(v)}\|f_1\|\|f_2\|.
		\]
		\end{enumerate}
 	\end{corollary}
	\begin{proof}
		(ii), (iii), and (iv) are equivalent by definition.
		The equivalence of (iv), (v), (vi), (vii), and (viii) is \cite[Corollary~1.3]{LWW}.
		Clearly, (ii) implies (i).
		To see (i) implies (iii) we apply Theorem~\ref{thm:muregularimplieslimitcone}
		to obtain $\langle\mu_\Gamma,\alpha\rangle=0$ for every $\alpha \in \Pi$.
		Hence, $\mu_\Gamma=0$.
	\end{proof}
	
	Since $\mathcal{L}_\Gamma' \subseteq \mathcal{L}_\Gamma$ we obtain Theorem~\ref{mainintro} from this corollary.
\begin{remark}
	We can indeed have $\mathcal{L}_\Gamma' = \emptyset$ but $\mathcal{L}_\Gamma \supsetneq \{0\}$.
		Let for example $G=G_1\times G_2$ and $\Gamma=\Gamma_1\times \Gamma_2$
		with $G_i$ of real rank one.
		As soon as $\Gamma_i$ are both infinite, $\mathcal{L}_\Gamma = \mathfrak{a}_+$
		and we cannot infer any regularity of $\mu_\Gamma$ from $\mathcal{L}_\Gamma$.
		However, if $\delta_{\Gamma_1} < \rho_1$
		and $\delta_{\Gamma_2} > \rho_2$, then
		$v_{\alpha_1} \notin \mathcal{L}_\Gamma'$
		where $\alpha_1$ is the unique reduced restricted root of $G_1$.
		Hence, $\langle\mu_\Gamma,\alpha_1\rangle = 0$ and therefore $\mu_\Gamma \in \R_{\geq 0} \alpha_2$.
		If $\delta_{\Gamma_1} < \rho_1$
		and $\delta_{\Gamma_2} < \rho_2$, then
		$\mathcal{L}_\Gamma' = \emptyset$.
	\end{remark}

Let us now formulate the consequences of Theorem~\ref{thm:muregularimplieslimitcone}
in analogy to the results of Subsection~\ref{ss:furtherresults} in the general (non Zariski dense) case.
	
 	\begin{corollary}
		\label{thm:onewallavoided}
		Let $\alpha\in \Pi$.
		\begin{enumerate}
			\item If $C_\alpha\cap \mathcal{L}'_\Gamma = \emptyset$,
				then
				\[
					\langle \mu_\Gamma,\alpha\rangle = \langle \mu_\Gamma,\iota\alpha\rangle = 0.
				\]
			\item If $C_\beta \cap \mathcal{L}_\Gamma'=\emptyset$
				for all $\beta\in \Pi \setminus\{\alpha,\iota\alpha\}$,
				then
				\[
				\mu_\Gamma = \max(0, \delta'_{\omega_\alpha + \iota\omega_\alpha}) (\omega_\alpha + \iota\omega_\alpha).
			\]
		\end{enumerate}
 	\end{corollary}

	We note that from Corollary~\ref{cor:twowallsavoided} we get the bound $\psi_\Gamma \leq \rho$
	by only considering the position of the modified limit cone $\mathcal{L}'_\Gamma$.
	Let us assume $C_\beta \cap \mathcal{L}'_\Gamma =\emptyset$ for all $\beta\in \Pi\setminus \{\alpha,\iota\alpha\}$
	for one simple root $\alpha$ as in Corollary~\ref{thm:onewallavoided} (ii).
	Then we can also get a bound of the growth indicator function.
	Indeed, we have $\psi'_\Gamma \leq \delta'_{\omega_\alpha + \iota \omega_\alpha}
		(\omega_\alpha +\iota \omega_\alpha)$.
		If $G$ is simple of higher rank, by \eqref{eq:propTpsi},
		$\psi'_\Gamma \leq \rho -\Theta$ for $\Theta$ the half sum of roots in
		an strongly orthogonal system.
		From this we can also get a bound on $\psi_\Gamma$
		without any assumption on the modified critical exponents $\delta'_\mu$.
		Indeed,
		since $\mu_\Gamma (v)\leq (\rho-\Theta)(v)$
		for all $v\in \mathfrak{a}_+$
		and $\mu_\Gamma = \delta'_{\omega_\alpha +\iota \omega_\alpha}
		(\omega_\alpha +\iota \omega_\alpha)$,
		\begin{equation}
			\label{eq:propTboundondelta}
			\delta'_{\omega_\alpha +\iota \omega_\alpha} 
			\leq \min_{v\in \mathfrak{a}_+} \frac{(\rho-\Theta)(v)}{(\omega_\alpha +\iota \omega_\alpha)(v)}.
		\end{equation}
		Thus we have proven the following proposition.
		\begin{proposition}
			\label{prop:boundwallavoided}
			If $G$ is simple of higher rank, $\Gamma$ is not a lattice,
			and $C_\beta \cap\mathcal{L}'_\Gamma =\emptyset$
			for all $\beta\in \Pi\setminus\{\alpha,\iota\alpha\}$,
			then
\[
	\psi_\Gamma \leq \rho + \left(\min_{v\in \mathfrak{a}_+} \frac{(\rho-\Theta)(v)}{(\omega_\alpha +\iota \omega_\alpha)(v)}\right) (\omega_\alpha +\iota \omega_\alpha).
\]
		\end{proposition}
		This bound is clearly not sharp but can improve the bound
		\eqref{eq:propTpsi}
		as the following example shows.
		
\begin{figure}[ht]
			\centering
			\label{fig:deltasB2}

			\includegraphics[width=\textwidth, trim = 1cm 8mm 3mm 5mm,clip]{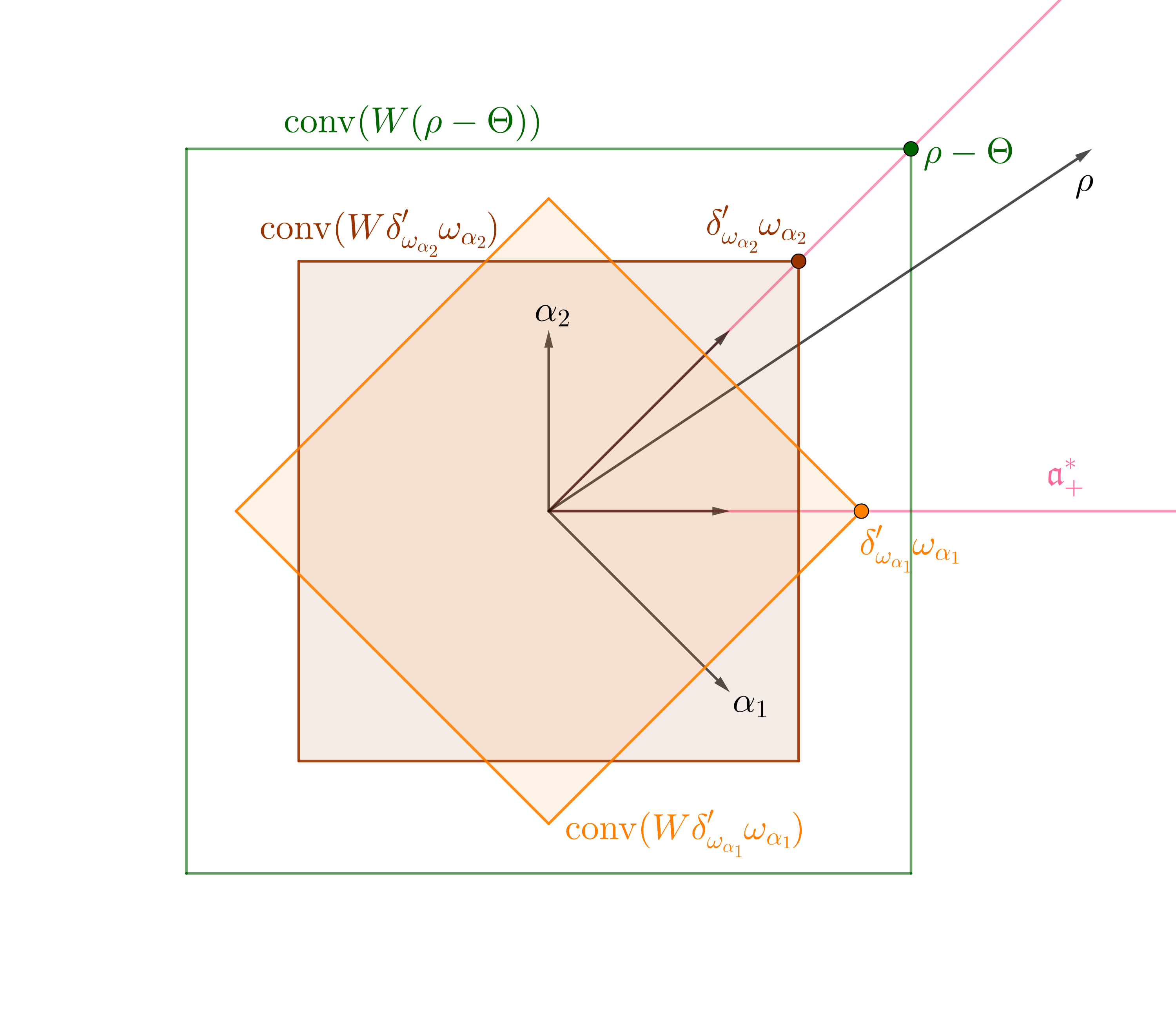}

			\caption{Visualisation of Proposition~\ref{prop:boundwallavoided}
				for $G=\mathrm{SO}_0(2,n)$.
				$\mu_\Gamma$ must be contained
				in $\conv{(\rho-\Theta)}$ (green).
				If $\mu_\Gamma \in \R_{\geq 0} \omega_{\alpha_2}$,
				then the maximal $\delta'_{\omega_{\alpha_2}}$ 
				is so that $\mu_\Gamma = \delta'_{\omega_{\alpha_2}} \omega_{\alpha_2}=\rho-\Theta$.
				There is no $\Gamma$-independent improvement since
				$\conv{(\rho-\Theta)}$ and $\conv{\delta'_{\omega_{\alpha_2}} \omega_{\alpha_2}}$
				(brown) coincide in this extremal case.
				Contrarily, if $\mu_\Gamma\in \R_{\geq 0}\omega_{\alpha_1}$,
				then $\conv{\delta'_{\omega_{\alpha_1}} \omega_{\alpha_1}}$ (orange)
				is always smaller than $\conv{(\rho-\Theta)}$
				as $\mu_\Gamma= \delta'_{\omega_{\alpha_1}} \omega_{\alpha_1}$ is contained in the latter.
			}
		\end{figure}

		\begin{example}
			Let $G=\mathrm{SO}_0(2,n)$, $n\geq 3$,
			and $\Gamma$ not a lattice.
			Then $\Sigma$ is of type $B_2$,
			i.e.~
			\[
				\Pi = \{\alpha_1 = (1,-1), \alpha_2 = (0,1)\} \subseteq 
				\mathfrak{a}^\ast
			\]
			acting on $\mathfrak{a} \simeq \R^2$.
			Moreover,
			\[
				\Sigma^+ =\{\alpha_1,\alpha_2,\alpha_3 = \alpha_1+\alpha_2 = (1,0), \alpha_4 = \alpha_1 + 2\alpha_2 = (1,1)\}
			\]
			and $\mathfrak{a}_+=\{v\in \R^2\colon v_1\geq v_2 \geq 0\}$.
		Furthermore, $\iota$ is trivial,
		\[
			\rho = \left(\frac n 2, \frac{n-2}2\right),
		\]
		and $\Theta = (1,0)$
		so that
		\[
			\rho - \Theta = \frac{n-2}2 (1,1).
		\]
		The fundamental weights are $\omega_{\alpha_1} = (1,0)=\alpha_3$
		and $\omega_{\alpha_2} = \frac 12 (1,1) = \frac 12 \alpha_4$.
		From \eqref{eq:propTboundondelta} we get
		\[
			\delta'_{\omega_{\alpha_1}} \leq 
			\min_{v\in \mathfrak{a}_+}
			\frac{(\rho-\Theta)(v)}{\omega_{\alpha_1}(v)}
			=
			\frac{n-2}2 \min_{v\in \mathfrak{a}_+}
			\frac{v_1+v_2}{v_1}= \frac{n-2}{2}
		\]
		and similarly
		\[
			\delta'_{2\omega_{\alpha_2}} \leq 
			\min_{v\in \mathfrak{a}_+}
			\frac{(\rho-\Theta)(v)}{2\omega_{\alpha_2}(v)}
			=
			\frac{n-2}2
			\min_{v\in \mathfrak{a}_+}
			\frac{v_1+v_2}{v_1+v_2}= \frac{n-2} 2.
		\]

		If $\mathcal{L}'_\Gamma \cap \ker \alpha_2 = \emptyset$,
		then Proposition~\ref{prop:boundwallavoided}
		yields
		\[
			\psi_\Gamma(v) \leq \rho(v) +
			\delta'_{2\omega_{\alpha_2}} 2\omega_{\alpha_2}(v) \leq 
			\rho(v) + \frac{n-2} 2 (v_1 + v_2)
			= (2\rho-\Theta)(v).
		\]
		Hence, this gives no improvement compared to \eqref{eq:propTpsi}
		as $\rho -\Theta$ and $\omega_{\alpha_2}$ are multiples of each other.

		In contrast, if $\mathcal{L}'_\Gamma \cap \ker \alpha_1 = \emptyset$,
		then $\psi_\Gamma \leq \rho + \frac{n-2}{2} \omega_{\alpha_1}$,
		i.e.~
		\[
			\psi_\Gamma(v) \leq (n-1) v_1 + \frac{n-2}{2}v_2.
		\]
	This improves \eqref{eq:propTpsi}
	which reads $\psi_\Gamma(v) \leq (n-1)v_1 + (n-2)v_2$.
	We refer to Figure~\ref{fig:deltasB2} for a visualization in terms of the convex hulls.
		\end{example}
		
	\section{The case of Zariski dense subgroups}
	In this section we apply the statements to the case
	where $\Gamma$ is Zariski dense.
	In this case $\psi_\Gamma$ is concave \cite{Qui02} and consequently
	$\mathcal{L}_\Gamma$ and $\mathcal{L}'_\Gamma$ are convex.
	The concavity is in fact the only ingredient that we will use in this section.
	We first note that
	\begin{equation}
		\label{eq:convcone}
		C_\alpha \cap \mathcal{L}'_\Gamma = \emptyset
		\iff
v_\alpha + v_{\iota\alpha}\notin \mathcal{L}'_\Gamma.
	\end{equation}
	Indeed, if $v\in C_\alpha \cap \mathcal{L}'_\Gamma$
	then we also have $\iota v\in C_\alpha \cap \mathcal{L}'_\Gamma$
	as $C_\alpha$ and $\psi_\Gamma$ are invariant under $\iota$.
	Thus $v + \iota v\in \mathcal{L}'_\Gamma$ by convexity.
	Since $v\in C_\alpha$, $v +\iota v$ is a multiple of $v_\alpha + v_{\iota\alpha}$
	which proves the equivalence \eqref{eq:convcone}.

	We now obtain Theorem~\ref{thm:oneweight} from \eqref{eq:convcone}, Corollary~\ref{thm:onewallavoided} (i),
and Proposition~\ref{prop:mu_is_biggest_growth}.
We then obtain Corollary~\ref{thm:onewallintro}
from the observation
\[
	\bigcap_{\beta\in \Pi\setminus\{\alpha,\iota\alpha\}} \ker \beta =  C_\alpha
\]
and
\[
	C_\alpha \cap \mathfrak{a}^{\mathrm{Her}} = \R_{\geq 0} (v_\alpha + v_{\iota\alpha}).
\]
It remains to prove Theorem~\ref{thm:psilinearintro}.
Here we use \eqref{eq:ofmainthm} to determine the growth indicator function $\psi_\Gamma$.

	\begin{proof}[Proof of Theorem~\ref{thm:psilinearintro}]
		If $\alpha\in I$, then $\langle\alpha,\mu_\Gamma\rangle >0$ and by Theorem~\ref{thm:muregularimplieslimitcone}
		there exists $v\in \mathcal{L}'_\Gamma \cap C_\alpha$
		such that $\psi'_\Gamma(v)=\mu_\Gamma(v)$.
		Using the concavity of $\psi_\Gamma'$ we find that
		\[
			\psi'_\Gamma(v+\iota v)\geq \psi'_\Gamma(v) + \psi'_\Gamma(\iota v)=2 \psi'_\Gamma(v)=2 \mu_\Gamma(v).
		\]
		On the other hand,
		\[
			\psi_\Gamma'(v+\iota v) \leq \mu_\Gamma(v+\iota v)
			=2\mu_\Gamma(v).
		\]
		As $v+\iota v \in \R_{\geq 0} (v_\alpha + v_{\iota\alpha})$ as before,
		we have $\psi'_\Gamma = \mu_\Gamma$
		on the ray $\R_{\geq 0}( v_\alpha + v_{\iota\alpha})$.

		The concavity implies that $\{\psi_\Gamma' \geq \mu_\Gamma\}$
		is a convex set $C$ containing $\R_{\geq 0} (v_\alpha + v_{\iota\alpha})$.
		Since $\psi'_\Gamma \leq \mu_\Gamma$ on $\mathfrak{a}_+$,
		it remains to see that $C\supseteq \mathfrak{a}_I \cap \mathfrak{a}^{\mathrm{Her}}_+$.
		Any $v\in \mathfrak{a}_+$ can be written as $\sum_{\alpha\in \Pi} c_\alpha v_\alpha$
		with $c_\alpha \geq 0$.
		The property $\iota v = v$ implies $c_\alpha = c_{\iota\alpha}$
		and $v\in \mathfrak{a}_I$ implies $c_\alpha = 0$ for $\alpha \notin I$.
		Hence, $v$ is a positive linear combination of $v_\alpha +v_{\iota\alpha}$ for $\alpha\in I$,
		and therefore contained in $C$.
		This completes the proof.
	\end{proof}

 	\begin{corollary}
		If $\Gamma$ is Zariski dense
		and $\mu_\Gamma$ is regular,
		then
		\[
			\psi_\Gamma=\mu_\Gamma +\rho
			\quad
		\text{on}
		\quad
			\mathfrak{a}_+^{\mathrm{Her}}
		\]
		and
		\[
			\mathfrak{a}_+^{\mathrm{Her}} \subseteq \mathcal{L}'_\Gamma \cup \{0\}.
		\]
 	\end{corollary}

	Let us also mention the special case where $\iota$ is trivial.
	This is for example the case if the root system $\Sigma$ of restricted roots
	is of type $B_n, C_n, D_n\: (n \text{ even}), E_7, E_8, F_4,$ or $G_2$.

 	\begin{corollary}
		If $\Gamma$ is Zariski dense,
		$\iota$ is trivial,
		and $\mu_\Gamma$ is regular,
		then
		\[
			\psi_\Gamma=\mu_\Gamma +\rho
			\quad
		\text{on}
		\quad
			\mathfrak{a}_+
		\]
		and
		\[
			\mathcal{L}'_\Gamma \cup \{0\}= \mathfrak{a}_+.
		\]
 	\end{corollary}
 \bibliographystyle{alpha}
 \bibliography{reference}

\begin{thebibliography}{BILW05}

\bibitem[AZ22]{AZ22}
J.-P. Anker and H.-W. Zhang.
\newblock Bottom of the {$L^2$} spectrum of the {L}aplacian on locally symmetric spaces.
\newblock {\em Geom. Dedicata}, 216(1):Paper No. 3, 12, 2022.

\bibitem[Ben96]{Ben96}
Y.~Benoist.
\newblock Actions propres sur les espaces homog\`enes r\'{e}ductifs.
\newblock {\em Ann. of Math. (2)}, 144(2):315--347, 1996.

\bibitem[Ben97]{Ben97}
Y.~Benoist.
\newblock Asymptotic properties of linear groups.
\newblock {\em Geom. Funct. Anal.}, 7(1):1--47, 1997.

\bibitem[BILW05]{BILW}
M.~Burger, A.~Iozzi, F.~Labourie, and A.~Wienhard.
\newblock Maximal representations of surface groups: symplectic {Anosov} structures.
\newblock {\em Pure Appl. Math. Q.}, 1(3):543--590, 2005.

\bibitem[BK15]{BK1}
Y.~Benoist and T.~Kobayashi.
\newblock Tempered reductive homogeneous spaces.
\newblock {\em J. Eur. Math. Soc. (JEMS)}, 17(12):3015--3036, 2015.

\bibitem[BK21]{BK3}
Y.~Benoist and T.~Kobayashi.
\newblock Tempered homogeneous spaces {III}.
\newblock {\em J. Lie Theory}, 31(3):833--869, 2021.

\bibitem[BK22]{BK2}
Y.~Benoist and T.~Kobayashi.
\newblock Tempered homogeneous spaces {II}.
\newblock In {\em Dynamics, geometry, number theory---the impact of {M}argulis on modern mathematics}, pages 213--245. Univ. Chicago Press, 2022.

\bibitem[BK23]{BK4}
Y.~Benoist and T.~Kobayashi.
\newblock Tempered homogeneous spaces {IV}.
\newblock {\em J. Inst. Math. Jussieu}, 22(6):2879--2906, 2023.

\bibitem[BL25]{BenoistLiang}
Y.~Benoist and S.~Liang.
\newblock On the rate of exponential decay of coefficients on homogeneous spaces.
\newblock Preprint, {arXiv}:2507.02469 [math.{GR}] (2025), 2025.

\bibitem[BT22]{BT22}
J.-P. Burelle and N.~Treib.
\newblock Schottky presentations of positive representations.
\newblock {\em Math. Ann.}, 382(3-4):1705--1744, 2022.

\bibitem[CHH88]{CHH88}
M.~Cowling, U.~Haagerup, and R.~Howe.
\newblock Almost {$L^2$} matrix coefficients.
\newblock {\em J. Reine Angew. Math.}, 387:97--110, 1988.

\bibitem[Cor90]{Cor90}
K.~Corlette.
\newblock {Hausdorff dimensions of limit sets I.}
\newblock {\em {Invent. Math.}}, 102(3):521--542, 1990.

\bibitem[CZZ22]{CZZ22}
R.~Canary, T.~Zhang, and A.~Zimmer.
\newblock Cusped {H}itchin representations and {A}nosov representations of geometrically finite {F}uchsian groups.
\newblock {\em Adv. Math.}, 404:Paper No. 108439, 67, 2022.

\bibitem[DK22]{DK}
S.~Dey and M.~Kapovich.
\newblock Patterson-{S}ullivan theory for {Anosov} subgroups.
\newblock {\em Trans. Am. Math. Soc.}, 375(12):8687--8737, 2022.

\bibitem[DKO24]{DKO}
S.~Dey, D.~Kim, and H.~Oh.
\newblock Ahlfors regularity of {P}atterson-{S}ullivan measures of {A}nosov groups and applications.
\newblock {\em arXiv:2401.12398, to appear in Compositio Math}, 2024.

\bibitem[DO25]{DO}
S.~Dey and H.~Oh.
\newblock Deformations of {Anosov} subgroups: limit cones and growth indicators.
\newblock {\em J. Lond. Math. Soc., II. Ser.}, 112(3):35, 2025.
\newblock Id/No e70280.

\bibitem[ELO23]{ELO23}
S.~Edwards, M.~Lee, and H.~Oh.
\newblock Anosov groups: local mixing, counting and equidistribution.
\newblock {\em Geom. Topol.}, 27(2):513--573, 2023.

\bibitem[Els73a]{MR360472}
J.~Elstrodt.
\newblock Die {R}esolvente zum {E}igenwertproblem der automorphen {F}ormen in der hyperbolischen {E}bene. {I}.
\newblock {\em Math. Ann.}, 203:295--300, 1973.

\bibitem[Els73b]{MR360473}
J.~Elstrodt.
\newblock Die {R}esolvente zum {E}igenwertproblem der automorphen {F}ormen in der hyperbolischen {E}bene. {II}.
\newblock {\em Math. Z.}, 132:99--134, 1973.

\bibitem[Els74]{MR360474}
J.~Elstrodt.
\newblock Die {R}esolvente zum {E}igenwertproblem der automorphen {F}ormen in der hyperbolischen {E}bene. {III}.
\newblock {\em Math. Ann.}, 208:99--132, 1974.

\bibitem[EO23]{OhTemperedness}
S.~Edwards and H.~Oh.
\newblock Temperedness of {$L^2(\Gamma \setminus G)$} and positive eigenfunctions in higher rank.
\newblock {\em Comm. Amer. Math. Soc.}, 3:744--778, 2023.

\bibitem[FO25]{FO25}
M.~Fr\k{a}czyk and H.~Oh.
\newblock Zariski-dense non-tempered subgroups in higher rank of nearly optimal growth.
\newblock {\em to appear in J. Reine Angew. Math.}, 2025.

\bibitem[GMT23]{OMT23}
O.~Glorieux, D.~Monclair, and N.~Tholozan.
\newblock Hausdorff dimension of limit sets for projective {Anosov} representations.
\newblock {\em J. {\'E}c. Polytech., Math.}, 10:1157--1193, 2023.

\bibitem[Hel84]{gaga}
S.~Helgason.
\newblock {\em Groups and geometric analysis: integral geometry, invariant differential operators, and spherical functions}.
\newblock Pure and applied mathematics. Academic Press, 1984.

\bibitem[KMO24]{KMO24}
D.M. Kim, Y.N. Minsky, and H.~Oh.
\newblock Tent property of the growth indicator functions and applications.
\newblock {\em Geom. Dedicata}, 218(1):18, 2024.
\newblock Id/No 14.

\bibitem[KT24]{KT}
F.~Kassel and N.~Tholozan.
\newblock Sharpness of proper and cocompact actions on reductive homogeneous spaces.
\newblock Preprint, {arXiv}:2410.08179, 2024.

\bibitem[Leu04]{Leu04}
E.~Leuzinger.
\newblock Critical exponents of discrete groups and {$L^2$}-spectrum.
\newblock {\em Proc. Amer. Math. Soc.}, 132(3):919--927, 2004.

\bibitem[LO23]{OhDichotomy}
M.~Lee and H.~Oh.
\newblock Dichotomy and measures on limit sets of {A}nosov groups.
\newblock {\em Int. Math. Res. Not. IMRN}, 2024(7):5658--5688, 10 2023.

\bibitem[LWW25]{LWW}
C.~Lutsko, T.~Weich, and L.L. Wolf.
\newblock Polyhedral bounds on the joint spectrum and temperedness of locally symmetric spaces.
\newblock Preprint, {arXiv}:2402.02530, to appear in Duke Math J., 2025.

\bibitem[Oh02]{OhUniformPointwise}
H.~Oh.
\newblock Uniform pointwise bounds for matrix coefficients of unitary representations and applications to {K}azhdan constants.
\newblock {\em Duke Math. J.}, 113(1):133--192, 2002.

\bibitem[Pat76]{pattersonlimitset}
S.J. Patterson.
\newblock The limit set of a {F}uchsian group.
\newblock {\em Acta Math.}, 136(3-4):241--273, 1976.

\bibitem[PS17]{PS17}
R.~Potrie and A.~Sambarino.
\newblock Eigenvalues and entropy of a {Hitchin} representation.
\newblock {\em Invent. Math.}, 209(3):885--925, 2017.

\bibitem[PSW23]{PSW23}
M.B. Pozzetti, A.~Sambarino, and A.~Wienhard.
\newblock Anosov representations with {Lipschitz} limit set.
\newblock {\em Geom. Topol.}, 27(8):3303--3360, 2023.

\bibitem[Qui02]{Qui02}
J.-F. Quint.
\newblock Divergence exponentielle des sous-groupes discrets en rang sup{\'e}rieur.
\newblock {\em Comment. Math. Helv.}, 77:563--608, 2002.

\bibitem[Sam14]{Sam}
A.~Sambarino.
\newblock Hyperconvex representations and exponential growth.
\newblock {\em Ergodic Theory Dyn. Syst.}, 34(3):986--1010, 2014.

\bibitem[Web08]{Web08}
A.~Weber.
\newblock Heat kernel bounds, {P}oincar\'{e} series, and {$L^2$} spectrum for locally symmetric spaces.
\newblock {\em Bull. Aust. Math. Soc.}, 78(1):73--86, 2008.

\bibitem[WW23]{WW23b}
T.~Weich and L.~L. Wolf.
\newblock Absence of principal eigenvalues for higher rank locally symmetric spaces.
\newblock {\em Communications in Mathematical Physics}, 403(3):1275--1295, 2023.

\bibitem[WW24]{WW23}
T.~Weich and L.~L. Wolf.
\newblock Temperedness of locally symmetric spaces: The product case.
\newblock {\em Geom. Dedicata}, 218(76), 2024.

\bibitem[WZ24]{WZ23}
L.~L. Wolf and H.-W. Zhang.
\newblock {$L^2$-spectrum, growth indicator function and critical exponent on locally symmetric spaces}.
\newblock {\em Proc. Amer. Math. Soc.}, 152:5445--5453, 2024.

\end{thebibliography}

\end{document}